\documentclass[10pt]{amsart}
\usepackage[english]{babel}
\usepackage{amscd, amssymb, amsmath, amsthm, graphics}
\usepackage{amsmath,amsfonts,amsthm,amssymb}
\usepackage{latexsym,amsmath}
\usepackage{mathrsfs}
\usepackage[all]{xy}

\newtheorem{theorem}{Theorem}[section]
\newtheorem{proposition}[theorem]{Proposition}
\newtheorem{corollary}[theorem]{Corollary}
\newtheorem{lemma}[theorem]{Lemma}

\theoremstyle{definition}
\newtheorem{definition}[theorem]{Definition}
\newtheorem{notation}[theorem]{Notation}
\newtheorem{example}[theorem]{Example}
\newtheorem{remark}[theorem]{Remark}

\theoremstyle{remark}
\newtheorem{question}{Question}

\numberwithin{equation}{section}

\renewcommand{\t}{ \widetilde}
\renewcommand{\hat}{ \widehat}

\newcommand{\Z}{\mathbb Z}
\newcommand{\Q}{\mathbb Q}
\newcommand{\R}{\mathbb R}

\newcommand{\N}{\mathbb N}
\newcommand{\C}{\mathbb C}

\newcommand{\tr}{{\rm tr}}

\newcommand{\ud}{{\mathrm{d}}}

\renewcommand{\o}{\overline}

\newcommand{\co}{\colon\thinspace}
\renewcommand{\epsilon}{\varepsilon}
\renewcommand{\c}{\mathcal}

\usepackage{times}

\begin{document}

\title{Volume of representations and mapping degree}  

\author{Pierre Derbez}
\address{Aix Marseille Univ, CNRS, Centrale Marseille, I2M, Marseille, France}
\email{pderbez@gmail.com}
\author{Yi Liu}
\address{Beijing International Center of Mathematics, Peking University, Beijing 100871, China}
\email{liuyi@math.pku.edu.cn}
\author{Hongbin Sun}
\address{Department of Mathematics, UC Berkeley CA 94720, USA}
\email{hongbins@math.berkeley.edu}
\author{Shicheng Wang}
\address{Department of Mathematics, Peking University, Beijing 100871, China}
\email{wangsc@math.pku.edu.cn}

\thanks{
The author Y.~L.~is partially supported by the Recruitment Program of Global Youth Experts of China.
The author H.~S.~is partially supported by Grant No. DMS-1510383 of the National Science Foundation of the United States.
The author S.~W.~is partially supported by grant No.~11371034 of the National Natural Science Foundation of China.}


\date{\today}
\begin{abstract}
Given a connected real Lie group and a contractible homogeneous proper $G$--space $X$ 
furnished with a $G$--invariant volume form, a real valued volume can be assigned to 
any representation $\rho\colon \pi_1(M)\to G$ for any oriented closed smooth manifold $M$ of the same dimension as $X$.
Suppose that $G$ contains a closed and cocompact semisimple subgroup, it is shown in this paper that 
the set of volumes is finite for any given $M$. From a perspective of model geometries,
examples are investigated and applications with mapping degrees are discussed.
\end{abstract}
\maketitle

\section{Introduction}
Let $G$ be a connected real Lie group and $H$ be any maximal compact subgroup of $G$.
After fixing a $G$--invariant volume form on the homogeneous $G$--space $X=G/H$,
a volume $\mathrm{vol}_G(M,\rho)\in \R$ can be assigned to any oriented closed smooth manifold $M$
of the same dimension as $X$ with respect to a representation $\rho\co\pi_1(M)\to G$.
The value 
$$\mathrm{V}(M,G)=\sup_\rho\,|\mathrm{vol}_G(M,\rho)|$$
in $[0,+\infty]$ is called the \emph{$G$--representation volume} of $M$.
This invariant has been introduced and studied by R.~Brooks and W.~Goldman 
\cite{BG1,BG2,Go} as 
a geometrical analogue of the celebrated simplicial volume
of orientable closed manifolds due to M.~Gromov \cite{Gromov,Th1}. 
During the past few years, much has been known about the $(\t{\mathrm{SL}}_2(\R)\times_{\Z}\R)$--representation volume (the Seifert volume)
and the $\mathrm{PSL}(2,\C)$--representation volume (the hyperbolic volume) for $3$--manifolds and their finite covers
\cite{DW1,DW2,DLW,DSW,DLSW}.
Those invariants have  demonstrated to be useful in studying nonzero degree maps between $3$--manifolds,
especially when the simplicial volume vanishes.
Now it seems a right time to consider higher dimensions and give a comprehensive treatment
of the theory under reasonably general hypotheses.

In this paper, we are primarily interested 
in the representation volume for Lie groups that contain closed and cocompact
semisimple subgroups. 
We formulate all the results in the smooth category to avoid technicalities.
Our main conclusion is that in this case, the representation volume
is a finite and nontrivial invariant 
for closed orientable smooth manifolds of the applicable dimension: 

\begin{theorem}\label{main-volume}
	Suppose that $G$ is a connected real Lie group
	which contains a closed cocompact connected semisimple subgroup. 
	Let $X=G/H$ be a homogeneous space furnished with a $G$--invariant volume form,
	where $H$ is a maximal compact subgroup of $G$.
	Then for any oriented closed smooth manifold $M$ of the same dimension 
	as $X$,	the volume function
	$$\mathrm{vol}_G\co \mathcal{R}(\pi_1(M),G)\to\R$$
	takes only finitely many values on the space of representations $\mathcal{R}(\pi_1(M),G)$.
	Moreover, there exists some aspherical $M$ for which
	$\mathrm{vol}_G$ is not constantly zero.
\end{theorem}

%

The assumptions here are natural from the perspective of maximal geometries with compact models.
According to W.~Thurston \cite{Th1}, a \emph{geometry} of dimension $n$ is a pair $(X,G)$
where $X$ is a simply-connected smooth manifold of dimension $n$ and 
$G$ is a real Lie group that acts transitively and effectively on $X$ by diffeomorphisms 
with compact isotropy at every point. It is a \emph{maximal geometry} if $G$ cannot be extended to
any larger groups, and it is a (compact) \emph{model geometry} 
if there is a closed $n$--manifold locally modeled on $(X,G)$. 
Two geometries are considered to be equivalent if there is an equivariant diffeomorphism
between the spaces with respect to the structure groups.
For classifications of geometries of dimension at most 5, 
we refer the reader to \cite{Scott,Hillman,Geng}.
Furnishing $X$ with a $G$--invariant volume form $\omega_X$,
we can speak of the representation volume of closed orientable smooth $n$--manifolds
with respect to the identity component $G^\circ$ of $G$.

\begin{corollary}\label{main-volume-corollary}
	Let $(X,G')$ be an $n$--dimensional geometry where $X$ is contractible 
	and $G'$ is semisimple.
	For any geometry $(X,G)$ that extends $(X,G')$ and 
	any $G$--invariant volume form $\omega_X$,
	the $G$--representation volume is a homotopy-type invariant 
	of closed orientable smooth $n$--manifolds,
	which takes values in $[0,+\infty)$ and does not always vanish.
\end{corollary}

The hyperbolic plane $(\mathbf{H}^2,\mathrm{Iso}(\mathbf{H}^2))$ is a well-known nontrivial example in dimension $2$. 
With respect to the hyperbolic metric,
it yields the $\mathrm{PSL}(2,\R)$--representation volume
which equals $-2\pi\chi(\Sigma)$ for any orientable closed surface $\Sigma$ of genus at least $2$, or $0$ otherwise.
Higher dimensional hyperbolic spaces $(\mathbf{H}^n,\mathrm{Iso}(\mathbf{H}^n))$ gives rise to
the $\mathrm{SO}^+(n,1)$--representation volume. 
It recovers the usual hyperbolic volume for orientable closed hyperbolic $n$--manifolds,
whereas for other $n$--manifolds the volume may be trivial or not.
Other symmetric spaces of noncompact types provide geometric models that
yield the representation volume for semisimple Lie groups without compact normal subgroups.
In dimension $3$, the $(\t{\mathrm{SL}}_2(\R)\times_{\Z}\R)$--representation volume arises 
from the Seifert-space geometry $(\t{\mathrm{SL}}_2(\R),\mathrm{Iso}(\t{\mathrm{SL}}_2(\R))$.
It can be greater than the $\t{\mathrm{SL}}_2(\R)$--representation volume in general.

In dimension $5$, there is an another interesting family which are non-symmetric spaces, 
namely, the maximal model geometries 
$(\t{\mathrm{SL}}_2(\R)\times_\alpha \t{\mathrm{SL}}_2(\R),\mathrm{Iso}(\t{\mathrm{SL}}_2(\R)\times_\alpha \t{\mathrm{SL}}_2(\R)))$
for $\alpha\in\Q^*$.
The identity component of the structure group is isomorphic to
$$(\t{\mathrm{SL}}_2(\R)\times_\Z\R)\times (\t{\mathrm{SL}}_2(\R)\times_\Z\R)/\{(t,\alpha t)\in\R\times \R\co t\in\R\},$$
which is $7$ dimensional.
If one regards the Seifert-space geometry $\t{\mathrm{SL}}_2(\R)$
as an affine real line bundle over the hyperbolic plane
via the fibration $\t{\mathrm{SO}}(2)\to \t{\mathrm{SL}}_2(\R)\to \mathbf{H}^2$,
the contractible homogeneous space $\t{\mathrm{SL}}_2(\R)\times_\alpha \t{\mathrm{SL}}_2(\R)$
can be observed as an affine real line bundle over $\mathbf{H}^2\times\mathbf{H}^2$,
where the fibers are modelled on the quotient of the affine real plane $\t{\mathrm{SO}}(2)\times\t{\mathrm{SO}}(2)$
by translations parallel to the $(1,\alpha)$ direction.
When $\alpha\in\R$ is irrational, one could also define $5$--dimensional geometries
with the identity component constructed by the same expression. 
The corresponding maximal analytic semisimple subgroup is no longer closed.
Nevertheless, we do not know if the irrational cases admit compact models anyways,
so there is a chance that they only yield vanishing representation volumes.

As a classical application, representation volumes can be used to study mapping degrees.
We say that a closed orientable smooth $n$--manifold $N$ has \emph{Property D} if 
for all closed orientable smooth manifold $M$ of the same dimension,
the set of mapping degrees
	$$\mathcal{D}(M,N)=\{|\mathrm{deg}(f)|\co f\in C^\infty(M,N)\}$$
is always finite.
In dimension 2, manifolds with Property D are those of the geometry $\mathbf{H}^2$,
by considering the simplicial volume.
In dimension 3, manifolds with Property D are those which contain non-geometric prime factors,
or prime factors of the geometry $\mathbf{H}^3$ or $\t{\mathrm{SL}}_2(\R)$,
by \cite[Corollary 1.7]{DLSW}.


By Corollary \ref{main-volume-corollary} 
and the domination inequality of representation volumes (Corollary \ref{comparison_of_volume}),
we immediately infer the following:

\begin{theorem}\label{modeled_Property_D}.
	If $(X,G)$ is a geometry where $X$ is contractible 
	and $G$ is semisimple, then every orientable closed manifold locally modeled on $(X,G)$
	has Property D.
\end{theorem}

Note that $X$ is a symmetric space of  non-compact type when $G$ has trivial center.
In that case, the result has been proved by Besson--Courtois--Gallot  \cite{BCG1} in rank one
and by Connell--Farb \cite{CF} in higher ranks (with a few exceptions), and it also 
follows from a stronger result that such manifolds have positive simplicial volume
by Lafont--Schmidt \cite{LS}. 
Both \cite{CF} and  \cite{LS}  extend the barycenter method of  \cite{BCG1}  from ranks one to higher ranks.
In other cases, Theorem \ref{modeled_Property_D} is actually talking about
a virtual (hyper-)torus bundle over a locally symmetric space, 
so the simplical volume has to vanish \cite[Section 3.1]{Gromov}.

%

In fact, except for the impossible dimensions $1$ and $2$ and the unknown dimension $4$,
there are always $n$--manifolds with Property D
which cannot be detected by the simplical volume.
This result has recently been established by Fauser and Loeh \cite[Theorem 1.3]{FL}.
In our terms, their examples are product manifolds $P\times Q$ such that $P$ and $Q$ are manifolds with Property D, 
and they show that $P\times Q$ must also have Property D.
If we take $P$ to be locally modeled on $\mathbf{H}^{n-3}$, 
and $Q$ to be locally modeled on $\t{\mathrm{SL}}_2(\R)$, then $P\times Q$ is locally modeled on $\mathbf{H}^{n-3}\times \t{\mathrm{SL}}_2(\R)$ which satisfies the assumptions of Theorem \ref{modeled_Property_D},
and we obtain an alternate proof of that result. 



We would like to make a remark about Theorem \ref{main-volume} 
before we close the introduction.
The prototype of Theorem \ref{main-volume} is of course the case of $(\t{\mathrm{SL}}_2(\R)\times_{\Z}\R)$--representations,
as established by Brooks--Goldman \cite{BG1,BG2}.
Their original proof is essentially reduced to the rigidity of certain characteristic class
associated with certain $\mathrm{PSL}(2,\R)$--representations, (precisely
the Godbillon--Vey class 
of the associated transversely--$\R P^1$ foliated circle bundles).
After their work, generalizations in various directions and alternate approaches 
are obtained for many semisimple Lie groups $G$.
Generalizations include
rigidity of characteristic classes for transversely $G$--homogeneous foliations,
and $G$--representation volumes for lattices, see \cite{Heitsch,BE}. 
One of the new approaches related to $3$--dimensional geometries identifies 
representation volumes for $\mathrm{PSL}(2,\C)$ and $\t{\mathrm{SL}}_2(\R)$
with the Cheeger--Chern--Simons invariants \cite{CnS, CgS}, 
see \cite{KK,Kh,DLW} and their references. 
This relation is first pointed out by 
Thurston for discrete and faithful representations in $\mathrm{PSL}(2,\C)$
of hyperbolic 3-manifold groups \cite{Th2}. 
Another new approach is due to Besson--Courtois--Gallot \cite{BCG2} 
based on the method they established in \cite{BCG1}.
The Lie groups considered there are the structure groups of 
rank-one symmetric spaces of non-compact type,
and a rigidity theorem for $\mathrm{Iso}_+(\mathbf{H}^n)$ is obtained.

Aiming at finding model geometries that naturally gives rise to representation volumes,
we have to extend our considerations to Lie groups that are not semisimple.
Indeed, a significant part of our argument provides more or less a classification of those groups
as assumed by Theorem \ref{main-volume}, (Proposition \ref{cocompact_coextension}).
Using that, we can reduce the proof of Theorem \ref{main-volume}
to a more concrete family, namely, the full central extension of real connected semisimple Lie groups 
with torsion-free center, (see Section \ref{Sec-centralExtension}).
The core argument for the central extension case is contained in the proof of Theorem \ref{rigidity}.
The proof of Theorem \ref{main-volume} is summarized in Section \ref{Sec-volumeFinitenessAndNontriviality}. 

The organization of this paper is reflected by the table of contents.
Efforts have been made to keep the exposition easily accessible to non-experts of Lie groups,
and the arguments as self-contained as possible.

\subsection*{Acknowledgement} 
We thank Jinpeng An for pointing us to a useful reference during the development of this paper.

\tableofcontents

\section{Volume of representations}\label{Sec-volumeOfRepresentations}
 
	In this preliminary section, we revisit this concept and
	scrutinize the construction in a natural setting.
	
	Let $G$ be any connected real Lie group. 
	Let $X$ be a contractible homogeneous proper $G$--space, namely, 
	a contractible smooth manifold equipped with 
	a transitive, proper action of $G$ by diffeomorphisms.
	Suppose that $\omega_X$ is a $G$--invariant volume form on $X$.
	Given any such triple $(G,X,\omega_X)$, 
	for any closed oriented smooth manifold $M$ of the same dimension as $X$,
	a real-valued function called \emph{volume} 
	can be defined over the $G$--representation variety $\c{R}(\pi_1(M),G)$ of $\pi_1(M)$.
	
	We confirm in this section that the volume function is continuous
	with respect to a natural topology of $\c{R}(\pi_1(M),G)$, 
	namely, the algebraic-convergence topology.
	Moreover, the volume function is essentially determined by $G$,
	up to a nonzero scalar factor proportional to the volume form.
	
	In the literature, 
	the $G$--space $X$ is usually taken to be $G/H$ 
	where $H$ is any maximal compact subgroup of $G$,
	so a $G$--invariant volume form $\omega_X$ can be given by any $G$--invariant Riemannian metric.
	We first point out the equivalence of that more concrete approach with our abstract setting.
	
	\begin{lemma}\label{concrete vs abstract}
		Let $G$ be any connected real Lie group
		and let $X$ be a contractible homogeneous proper $G$--space.
		Then 
		\begin{enumerate}
		\item	the point stabilizers of $X$ are maximal compact subgroups of $G$. 
		Hence $X$ is determined by $G$ up to a $G$--equivariant diffeomorphism,
		\item there exists a $G$--invariant volume form $\omega_X$ on $X$.
		Moreover, it is unique up to a signed rescaling.
		\end{enumerate}
	\end{lemma}
	
	\begin{proof}
		Recall that the action of a topological group $G$ on 
		a topological space $X$ is said to be \emph{proper} if the mapping $G\times X\to X\times X$ given by $(g,x)\mapsto(x,g\cdot x)$
		is proper, namely, such that the preimage of any compact subset is compact.
		Since the $G$--action on $X$ is assumed to be proper,
		point stabilizers $G_x$ for any $x\in X$ are compact subgroups of $G$, and the homogeneity assumption
		implies that $X$ is isomorphic to $G/G_x$ as a $G$--space.
		As any maximal compact subgroup is a topologically deformation retract of $G$,
		$X$ is homotopy equivalent to a closed manifold $H/G_x$,
		for some maximal compact subgroup $H$ that contains $G_x$.
		However, the assumption that $X$ is contractible forces $H/G_x$ to be a point, so $G_x$
		is a maximal compact subgroup of the Lie group $G$.
		As maximal compact subgroups of $G$ are conjugate to each other,
		different $G$--spaces $X$ as assumed are $G$--equivariantly diffeomorphic to each other.
		This proves  statement (1).
		
		To prove  statement (2), we may now take $X$ to be $G/H$ as usual, where $H$ is any maximal compact subgroup of $G$.
		In this case, it is well-known that there exists a complete $G$--invariant Riemannian metric on $G/H$.
		For example, one may first take an arbitrary inner product $\langle\cdot,\cdot\rangle'_o$ on the tangent space $T_oX$ 
		at the point $o\in X$ which corresponds to the coset $H$.
		Then average the inner product over $H$ by defining
			$$\langle V,W\rangle_o=\int_{H} \langle h_*V, h_*W\rangle'_o\,\mathrm{d}\mu_H(h),$$
		for all $V,W\in T_oX$, where $\mu_H$ denotes the normalized Haar measure of $H$.
		After that, a claimed $G$--invariant Riemannian metric
		can be obtained by translating the $H$--invariant inner product at $T_oX$ 
		over $X$ by the action of $G$.
		The associated volume form $\omega_X$
		is hence $G$--invariant as claimed.
		To see the uniqueness of the $G$--invariant volume form up to rescaling,
		it suffices to observe that any such form is determined by its restriction at $o$,
		namely, $\omega_X|_o$, which lies in $\wedge^{\dim X}T^*_oX\cong\R$.
		As $\omega_X|_o$ must be nonzero, any pair of $G$--invariant volume forms over $X$
		differ only by a multiplication of some nonzero real scalar.
		This proves  statement (2).	
		Note that, however, $G$--invariant Riemannian metrics on $X$ are not necessarily proportional, 
		unless the isotropy representation of $H$	on $T_oX$ is irreducible.
	\end{proof}
	
	Let $(G,X,\omega_X)$ be any triple as assumed. For any oriented (connected) closed smooth manifold $M$
	of the same dimension as $X$, we regard $\pi_1(M)$ as the (discrete) deck transformation group 
	of the universal cover $\widetilde{M}$ of $M$.
	
	Denote by $\c{R}(\pi_1(M),G)$ the $G$--representation variety of $\pi_1(M)$.
	In this paper, for any finitely generated group $\pi$,
	we understand the $G$--representation variety $\c{R}(\pi,G)$ 
	as the set of all  homomorphisms of $\pi$ into $G$.
	This set is inherited with the algebraic-convergence topology,
	which can be characterized by the property that a sequence 
	$\{\rho_n\}_{n\in\mathbb{N}}$
	converges to $\rho_\infty$ in $\c{R}(\pi,G)$
	if and only if $\{\rho_n(g)\}_{n\in\mathbb{N}}$
	converges to $\rho_\infty(g)$ in $G$ for every $g\in\pi$.
	
	The volume function 
		$${\rm vol}_{G,X,\omega_X}\colon \c{R}(\pi_1(M),G)\to\R$$
	can be defined as follows.
	For any representation
		$$\rho\colon \pi_1(M)\to G,$$
	denote by 
		$$M\times_\rho X$$
	the associated $G$--flat $X$--bundle space
	over $M$. It is the quotient of the product space $M\times X$
	by the diagonal freely discontinuous action of $\pi_1(X)$, namely,
	$\sigma\cdot(m,x)=(\sigma\cdot m,\rho(\sigma)\cdot x)$ for all $\sigma\in\pi_1(X)$.
	The bundle projection $M\times_\rho X\to M$ 
	is induced by the projection of $\widetilde{M}\times X$ onto the first factor;
	the flat structure is the induced horizontal foliation 
	on $M\times_\rho X$ with holonomy in $G$.
	The pull-back form of $\omega_X$ on $\widetilde{M}\times X$, via the projection onto the second factor,
	is invariant under the diagonal action of $\pi_1(M)$, 
	so there is a naturally induced form $\omega^{\rho}_X$ on $M\times_\rho X$.
	Take $s\colon M\to M\times_\rho X$ to be any $C^1$--differentiable section of the bundle.
	
	We define the \emph{volume} of $(M,\rho)$ with respect to 
	the triple $(G,X,\omega_X)$ to be
	\begin{equation}\label{vol_section}
		{\rm vol}_{G,X,\omega_X}(\rho)=\int_{M} s^*\omega^{\rho}_X.
	\end{equation}
		
	There is a natural way to translate between sections and developing maps. 
	Specifically, every section $s\colon M\to M\times_\rho X$
	induces unique a lift $\tilde{s}\colon \t{M}\to \t{M}\times X$, after choosing a base point lift of $M$ into $\widetilde{M}$.
	Then the projection of $\tilde{s}$ to the second factor
	gives rise to an equivariant map with respect to $\rho\colon \pi_1(M)\to G$,
		$$D_\rho\colon \widetilde{M}\to X,$$
	which is called a \emph{developing map}. Conversely, the graph of any developing map $D_\rho$ in $M\times X$
	is invariant under the action of $\pi_1(M)$, so $\mathrm{Id}_M\times D_\rho$ induces a section $s$. 
	In terms of a developing map $D_\rho$, we may therefore write
	\begin{equation}\label{vol_developing}
		{\rm vol}_{G,H,\omega_X}(\rho)=\int_{\c{F}} D^*_\rho\omega_X,
	\end{equation}
	where $\c{F}\subset \widetilde{M}$ is any fundamental domain of $\pi_1(M)$.	
		
	The dependence of the definition on our auxiliary choices can be summarized by the following:
	
	\begin{lemma}\label{choices}
		With the notations and assumptions above, the following statements are true:
		\begin{enumerate}
			\item A $C^1$--differentiable section $s$ of $M\times_\rho X$ exists and any such section
			yields one and the same value 
			$${\rm vol}_{G,X,\omega_X}(\rho)\in \R.$$
			\item For any two $G$--invariant volume forms $\omega_X$ and $\omega'_X$,
			the associated volume functions agree up to a nonzero scalar factor.
			Namely, there exists $\lambda\in \R^*$ such that
			$${\rm vol}_{G,X,\omega_X}(\rho)=\lambda\,{\rm vol}_{G,X,\omega'_{X}}(\rho)$$
			for all $M$ and all $\rho\in\c{R}(\pi_1(M),G)$.
		\end{enumerate}
	\end{lemma}
	
	\begin{proof}
		Since $X$ is contractible, the bundle $M\times_\rho X$ admits a continuous section
		and any two sections $s,s'\colon M\to M\times_\rho X$ are homotopic.
		Using standard approximation techniques in differential topology,
		the sections and homotopy can be approximated by $C^r$--differentiable ones,
		for arbitrary $1\leq r\leq \infty$,
		(see \cite[Chapter 2, Section 2, Exercise 3]{Hirsch}).
		The induced form $\omega^\rho_X$ is closed as is $\omega_X$.
		So the integral over $M$ against the pull-backs of $\omega^\rho_X$ via $s$ and $s'$ are equal.
		This proves  statement (1).
		
		Statement (2) follows immediately from Lemma \ref{concrete vs abstract} (2).
	\end{proof}
	
	The following proposition shows that the volume function is well behaved.
	The continuity is essentially a consequence	of the properness assumption about the $G$--action on $X$.
	
	\begin{proposition}\label{volume_continuity}
		The function ${\rm vol}_{G,X,\omega_X}\co \c{R}(\pi_1(M),G)\to \R$ is conjugation-invariant 
		and continuous.
	\end{proposition}
	
	\begin{proof}
		Given a developing map $D_\rho\co \t{M}\to X$ for a representation $\rho\in \c{R}(\pi_1(M),G)$, 
		it is straightforward to check that	$g\circ D_\rho$ 
		is a developing map for any conjugation $\tau_{g}(\rho)\in \c{R}(\pi_1(M),G)$,
		defined by $\tau_g(\rho)(\sigma)=g\rho(\sigma) g^{-1}$.
		The volume function is invariant under conjutation of representations by the $G$--invariance
		of the volume form $\omega_X$.
		
		%
		It remains to show the continuity.
		Fix a finite generating set of $\pi_1(M)$ induced by a finite cellular decomposition of $M$, and denote by $\mathcal{F}$
		a fundamental domain which is a connected union of lifted cells in $\widetilde{M}$.
		By the properness of the $G$--action on $X$, we can fix a $G$--invariant Riemannian metric of $X$, 
		as argued in Lemma \ref{concrete vs abstract}. Fix a Riemannian metric $g_M$ of $M$. 
		
		Suppose that $\{\rho_j\}_{j\in\N}$ is any sequence of representations
		which converges to a representation $\rho_\infty$.
		Take a $C^1$--differentiable developing map $D_\infty\colon \widetilde{M}\to X$ of  $\rho_\infty$.
		Then the algebraic convergence allows us to construct a sequence of $C^1$--differentiable developing maps
		$D_n\colon \widetilde{M}\to X$ which converge to $D_\infty$ uniformly over $\mathcal{F}$ 
		in the $C^1$--norm.
		For example, a sequence of $C^1$--uniformly convergent 
		developing maps $\{D_j\}_{j\in\N}$ can be built first over a neighborhood of the $0$--skeleton of $\mathcal{F}$,
		and then extended $\rho_j$--equivariantly over a neighborhood of the $0$--skeleton of $\widetilde{M}$, 
		and inductively over neighborhood of higher dimensional	skeleta of $\widetilde{M}$.
		Moreover, as $D_\infty$ is $C^1$--differentiable, we can require that in each step of the induction,
		the maps $D_j$ are constructed to be $C^1$--uniformly convergent to the $D_\infty$ 
		on the considered neighborhood of the skeleton.
		Provided with such a sequence of developing maps $\{D_j\}_{j\in\N}$, 
		given any constant $\epsilon>0$, 
		and for all sufficiently large $j$,
		we have the following estimation by the Stokes Formula:
			$$|{\rm vol}_{G,X,\omega_X}(\rho_j)-{\rm vol}_{G,X,\omega_X}(\rho_\infty)|<C\epsilon,$$
		where $C$ is a constant determined by the Lipschitz constant of $D_\infty$,
		and the area of $\partial\mathcal{F}$ with respect to $g_M$.
		As $\epsilon>0$ is arbitrary, we have
			$$\lim_{j\to\infty}{\rm vol}_{G,X,\omega_X}(\rho_j)={\rm vol}_{G,X,\omega_X}(\rho_\infty).$$
		This completes the proof.
	\end{proof}
		
	When one wishes to compute the volume for a $G$--representation, 
	it is important to keep in mind 
	that the exact value depends on a particular normalization of the volume measure.
	In certain cases, there are canonical ways to specify a normalization.
	For example, for volume of representations in $\mathrm{PSL}(2,\C)$,
	the normalization is customarily taken 
	as the volume form of hyperbolic geometry $\mathbf{H}^3$ (of constant curvature $-1$);
	for volume of representations in $\t{\mathrm{SL}}_2(\R)\times_{\Z}\R$, the normalization
	used by Brooks and Goldman \cite{BG1} is specified 
	in terms of the Godbillon--Vey invariant of the naturally associated transversely--$\R P^1$ foliated circle bundle.
	However, there is no natural way to specify a canonical normalization 
	to any connected Lie group in general.
	
	On the other hand, sometimes one may only be concerned with the ratio between volume values.
	In other words,	only those normalization-free qualities	of the volume function
	are considered to be interesting.
	%
	The discussion of this paper is of the latter flavor. 
	For this reason, we often adopt the following simpler and more conventional
	
	\begin{notation}\label{usualNotation}
		For any connected real Lie group $G$, we speak of volume of representations with a triple $(G,X,\omega_X)$ implicitly assumed as the following.
		The contractible proper homogeneous $G$--space $X$ is taken to be $G/H$ where $H$ is a fixed maximal compact subgroup of $G$.
		The $G$--invariant volume form is chosen and fixed. 		
		For any oriented closed smooth manifold $M$,
		and any representation $\rho\colon \pi_1(M)\to G$, we only retain the subscript $G$ and write
			$${\rm vol}_G(M,\rho)\in \R$$
		for the volume of representation ${\rm vol}_{G,X,\omega_X}(\rho)$.
	\end{notation}
	
	\begin{definition}
		Adopting Notation \ref{usualNotation}, we define the \emph{$G$--representation volume} of $M$ to be
			$${\rm V}(M,G)=\sup_{\rho\in\c{R}(\pi_1(M),G)}\,|{\rm vol}_G(M,\rho)|,$$
		which is valued in $[0,+\infty]$.
	\end{definition}

\section{Comparison of volumes}\label{Sec-comparisonOfVolumes}
	We consider the behavior of the volume of representations under variation of 
	the manifolds and the Lie groups.
	
	\begin{proposition}\label{change_of_volume}
		The following formulas hold for volume of representations.
		\begin{enumerate}
		\item Let $G$ be a connected real Lie group which acts properly and transitively on a contractible homogeneous $G$--space
		$X$ with a $G$--invariant volume form $\omega_X$. Let $M$ be any closed oriented smooth manifold of the same dimension as $X$.
		For any smooth map $f\co M'\to M$ where $M'$ is a closed oriented smooth manifold $M'$ of the same dimension as $X$,
		and any representation $\rho\co\pi_1(M)\to G$,
			$$\mathrm{vol}_{G}(M',f^*(\rho))=\mathrm{deg}(f)\cdot\mathrm{vol}_{G}(M,\rho),$$
		where $\mathrm{deg}(f)\in\Z$ denotes the signed mapping degree of $f$.
		\item Let $(G,X,\omega_X)$ and $M$ be the same as above.
		For any connected real Lie group $G'$ and any homomorphism $\phi\co G'\to G$
		with compact kernel and with closed and cocompact image,
		the induced action of $G'$ on $X$ is proper and transitive and preserves $\omega_X$.
		Moreover, for any representation $\rho'\co \pi_1(M)\to G'$,
			$$\mathrm{vol}_{G}(M,\phi_*(\rho'))=\mathrm{vol}_{G'}(M,\rho').$$
		\item Let $(G,X,\omega_X)$ be the same as above and let $I$ be a finite set.
		For all $i\in I$, let $M_i$ be a closed oriented smooth manifold of the same dimension  as $X$,
		and $\rho_i\co \pi_1(M_i)\to G$ be a representation.
		Then 
			$$\mathrm{vol}_G(\#_iM_i,\,\bigvee_i\rho_i)=\sum_i\mathrm{vol}_G(M_i,\rho_i).$$
		\item Let $I$ be a finite set and for all $i\in I$, let $(G_i,X_i,\omega_{X_i})$ be triples similar as above.
		Let $M_i$ be closed oriented smooth manifolds of the same  dimension  as $X_i$,
		and let $\rho_i\co \pi_1(M_i)\to G_i$ be representations. 
		Then with respect to the triple $(\prod_iG_i,\prod_iX_i,\prod_i\omega_{X_i})$,
			$$\mathrm{vol}_{\prod_iG_i}\left(\prod_iM_i,\prod_i\rho_i\right)=\prod_i\mathrm{vol}_{G_i}(M_i,\rho_i).$$
		\end{enumerate}
	\end{proposition}
	
	\begin{proof}
		To prove (1), 
		we observe that if 
		a developing map $D_\rho\colon \t{M}\to X$ defines the section $s\colon M\to M\times_{\rho}X$,
		then the composed map $D_\rho\circ \t{f}\colon \t{M'}\to\t{M}\to X$ is a developing map for $f^*(\rho)$,
		which defines a section	$s'\colon M'\to M'\times_{\rho\circ f_{\sharp}}X$.
		The form $(s')^*\omega^{\rho\circ f_\sharp}_X$ on $M'$ is induced by 
		the pull-back of $\omega_X$ to $\t{M}'$ via $D_\rho\circ\t{f}$,
		and it coincides with $(s\circ f)^*\omega^\rho_X$ by definition.		
		Turning the expression (\ref{vol_section}) into pairing between homology and cohomology classes, we have
			$$\mathrm{vol}_{G}(M',f^*(\rho))=\langle [M'], (s')^*[\omega^{\rho\circ f_\sharp}_X]\rangle
			=\langle [M'],(s\circ f)^*[\omega^\rho_X]\rangle=\langle f_*[M'], s^*[\omega^\rho_X]\rangle.$$
		As $f_*[M']=\mathrm{deg}(f)[M]$, we have
			$$\mathrm{vol}_{G}(M',f^*(\rho))=\mathrm{deg}(f)\cdot\langle [M],s^*[\omega^\rho_X]\rangle=\mathrm{deg}(f)\cdot\mathrm{vol}_{G}(M,\rho).$$
			
		To prove (2), denote by $K$ the kernel of $\phi$.
		Take $H$ to be a maximal compact subgroup of $G$, for example,
		the isotropy group of the action of $G$ on $X$. Since $K$ is assumed to be compact, 
		and $\phi(G')$ is assumed to be closed in $G$, $H'=\phi^{-1}(H)$ is compact in $G$.
		This means that $G'$ acts properly on $X$ via $\phi$.
		Moreover, the maximality of $H$ implies that $H'$ must be a maximal compact subgroup of $G$.
		In particular, if we naturally identify $X$ with the coset space $G/H$, 
		the subspace $\phi(G')H/H$ is an embedded closed contractible subspace of $X$. 
		Here $\phi(G')H$ denotes the double coset contained by $G$ 
		which consists of elements of the form $\phi(g')h$
		for all $g'\in G'$ and $h\in H$.
		By homotopy equivalences $G/\phi(G')\simeq G/\phi(H')\simeq H/\phi(H')$, 
		we can count the cohomological dimension of the compact spaces
		$H$, $\phi(H')$, and $G/\phi(G')$ by our assumption.
		Then we have $\dim G/\phi(G')=\dim H-\dim \phi(H')$. 
		This means $\phi(G')H/H$ equals $X$,
		or in other words, $G'$ acts transitively on $X$ via $\phi$.
		Of course the induced action of $G'$ on $X$ preserves the volume form $\omega_X$.
		
		Observe that for any representation $\rho'\colon \pi_1(M)\to G'$,
		any $\phi_*(\rho')$--equivariant developing map $\t{M}\to X$ is by definition 
		a $\rho'$--equivariant developing map for the induced action of $G'$, so
		$$\mathrm{vol}_{G}(M,\phi_*(\rho'))=\mathrm{vol}_{G'}(M,\rho').$$
		
		To prove (3), we take a developing map  which sends the spheres that the connected sum
		identifies along to points in $X$, then the claimed formula follows 
		immediately from the formula \ref{vol_developing}.
		
		To prove (4), we take a developing map of the form $\prod_iD_{\rho_i}\co \prod_i\t{M}\to \prod_i X_i$,
		and again we apply the formula \ref{vol_developing} to derive the claimed formula.
	\end{proof}
	
	\begin{corollary}\label{comparison_of_volume}
		The following comparisons hold for the representation volume of manifolds:
		\begin{enumerate}
		\item \emph{(Domination Inequality).} Given any $f\co M'\to M$ as of Proposition \ref{change_of_volume} (1),
			$$\mathrm{V}(M',G)\geq|\mathrm{deg}(f)|\cdot\mathrm{V}(M,G).$$
		\item \emph{(Induction Inequality).} Given any $\phi\co G'\to G$ as of Proposition \ref{change_of_volume} (2),
			$$\mathrm{V}(M,G')\leq \mathrm{V}(M,G).$$
		\item \emph{(Connected-Sum Inequality).} Given any $(G,X,\omega_X)$ and $M_i$ as of Proposition \ref{change_of_volume} (3),
			$$\mathrm{V}(\#_iM_i,\,G)\leq\sum_i\mathrm{V}(M_i,G).$$
		\item \emph{(Product Inequality).} Given any $(G_i,X_i,\omega_{X_i})$ and $M_i$ as of Proposition \ref{change_of_volume} (4),
			$$\mathrm{V}\left(\prod_iM_i,\prod_iG_i\right)\geq\max_{\sigma\in\mathfrak{S}(I)}\left\{\prod_i\mathrm{V}\left(M_i,G_{\sigma(i)}\right)\right\},$$
			where $\mathfrak{S}(I)$ is the permutation group of the finite set $I$ of indices,
			and the volume $\mathrm{V}(M_i,G_{\sigma(i)})$ is considered to be $0$ if 
			the dimension of $M_i$ mismatches that of $X_{\sigma(i)}$.
		\end{enumerate}
	\end{corollary}
	%
	
	None of the above inequalities are equalities
	in general.	In fact, the following examples show that each of them 
	can fail to be equal under certain circumstances.
	One can certainly find tons of such examples, 
	some of which may be trickier to verify.
	To avoid technicalities, 
	we only sketch the key points to verify our examples.
	The reader may safely skip this part for the first reading.
		
	\begin{example}\
	\begin{enumerate}
	\item Take $G$ to be either $\mathrm{PSL}(2,\C)$ or $\t{\mathrm{SL}}_2(\R)$.
	By \cite[Corollary 1.8]{DLW},
	there exists a closed oriented $3$--manifold $M$ with vanishing $\mathrm{V}(M,G)$,
	whereas $\mathrm{V}(M',G)>0$ holds for some finite cover $M'$ of $M$.
	So a strict domination inequality 
	$\mathrm{V}(M',G)>|\mathrm{deg}(f)|\,\mathrm{V}(M,G)$ holds for the covering map $f\co M'\to M$.
	\item Take $G$ to be $\t{\mathrm{SL}}_2(\R)$ and $G'$ to be $\t{\mathrm{SL}}_2(\R)\times_\Z\R$.
	Consider an oriented closed 3--manifold $M$ which is a Seifert fibered space with the symbol
	$(2,0;3/2)$. The base $2$--orbifold is a closed oriented surface of genus $2$ and with a cone point of order $2$,
	so it has Euler characteristic $\chi=-5/2$. The (orbifold) Euler number of the Seifert fibration 
	equals $e=3/2$. Then $\mathrm{V}(M,G')=4\pi^2\times (-5/2)^2/|3/2|=50\pi^2/3$ 
	by the formula of Seifert volume in this case, (see \cite[Section 6.2]{DLW}).
	For any representation $\rho\colon \pi_1(M)\to\t{\mathrm{SL}}_2(\R)$ of nonzero volume,
	it can be calculated by \cite[Proposition 6.3]{DLW} that $\rho$ must send the regular fiber
	$h$ to a generator $\mathrm{sh}(\pm1)$ of the center $\{\mathrm{sh}(n)\in \t{\mathrm{SL}}_2(\R)\colon n\in\Z\}$.
	Then $\mathrm{V}(M,G)= 4\pi^2\times 1^2\times(3/2)=6\pi^2$.
	So a strict induction inequality $V(M,G')<V(M,G)$ 
	holds for this case.
	\item	Take $G$ to be $\mathrm{PSL}(2,\C)$ and $\Gamma$ be a torsion-free uniform lattice of $G$.
	Denote by $M=\mathbf{H}^3/\Gamma$ be the closed hyperbolic $3$--manifold with the induced orientation.
	By Theorem \ref{main-volume}, the set of volumes for all representations of $\pi_1(M)\cong \Gamma$ in $G$ consists of 
	finitely many real values $v_1<v_2<\cdots<v_s$.
	If we require further that $M$ admits no orientation-reversing self-homeomorphism,
	it is implied by the volume rigidity that $v_s=\mathrm{Vol}_{\mathbf{H^3}}(M)$ and $|v_1|<v_s$,
	 \cite[Theorem 1.4]{BCG2}.
	Note that the set of volumes of the orientation-reversal $-M$ consists of $-v_s<\cdots<-v_2<-v_1$.
	Then the set of volumes for the connected sum $M\# (-M)$ consists of all the values $v_i-v_j$,
	and $\mathrm{V}(M\#(-M),G)=|v_s-v_1|<2v_s$.
	As $\mathrm{V}(M,G)=\mathrm{V}(-M,G)=v_s$,
	we obtain a strict connected-sum inequality 
	$\mathrm{V}(M\#(-M),G)<\mathrm{V}(M,G)+\mathrm{V}(-M,G)$
	in this case.
	\item Take $G$ to be $\t{\mathrm{SL}}_2(\R)\times\mathrm{PSL}(2,\C)$. 
	Let $S$ be the unit tangent bundle of a closed oriented hyperbolic surface,
  and $H$ be a closed oriented hyperbolic $3$--manifold.
	Let $M_S=S\times S$ and $M_H=H\times H$.
	
	We argue that $\mathrm{V}(M_S,G)$ must be zero. This means that
	$\mathrm{vol}_G(M_S,\rho)=0$ holds for every representation $\rho \co \pi_1(M_S)\to G$.
	Every $\rho$ is of the form $(\sigma,\eta)$ 
	where $\sigma\co\pi_1(M_S)\to \t{\mathrm{SL}}_2(\R)$,
	and $\eta\co \pi_1(M_S)\to\mathrm{PSL}(2,\C)$. 
	Via any section 
	$s=(u,v)\co M_S\to M_S\times_{(\sigma,\eta)}(\t{\mathrm{SL}}_2(\R)\times \mathbf{H}^3)$,
	the pull-back of the product volume form can be regarded as a cup product
	$[u^*\omega_{\t{\mathrm{SL}}_2(\R)}^\sigma]\smile[v^*\omega_{\mathbf{H}^3}^\eta]$
	in $H^6(M_S;\R)$,
	where $[u^*\omega_{\t{\mathrm{SL}}_2(\R)}^\sigma]$ and $[v^*\omega_{\mathbf{H}^3}^\eta]$
	are the $3$--dimensional factor forms lying in $H^3(M_S;\R)$.
	Note that	$H_3(M_S;\R)\cong \bigoplus_{i=0}^3(H_{3-i}(S;\R))\otimes H_i(S;\R))$.
	Moreover, the $(3,0)$--summand is generated by the $3$--cycle $S\times\mathrm{pt}$; the $(2,1)$--summand
	is generated by all the $3$--cycles of the form $T_\alpha\times\beta$ where $T_\alpha$ is any vertical torus of $S$
	over a non-separating simple closed curve of the base,	
	and $\beta$ is any loop of $S$ that projects to a non-separating simple closed curve of the base;
	the $(1,2)$--summand and the $(0,3)$--summand can be described similarly.
	The evaluation of $[u^*\omega_{\mathbf{H}^3}^\eta]$ on any of the above cycles 
	is nothing but the volume of the $\mathrm{PSL}(2,\C)$--representations 
	of their fundamental groups induced by $\eta$, so it is always trivial by a direct argument.
	It follows that $[u^*\omega_{\mathbf{H}^3}^\eta]=0$ in $H^3(M_S;\R)$, 
	so $[u^*\omega_{\t{\mathrm{SL}}_2(\R)}^\sigma]\smile[v^*\omega_{\mathbf{H}^3}^\eta]=0$
	in $H^6(M_S;\R)$. Therefore, $\mathrm{vol}_G(M_S,\rho)=0$ holds as claimed.
	
	However, it is obvious that $\mathrm{V}(M_S\times M_H,G\times G)>0$.
	Then we obtain a strict product inequality
	$\mathrm{V}(M_S\times M_H,G\times G)>\mathrm{V}(M_S,G)\times\mathrm{V}(M_H,G)$
	in this case.	
	\end{enumerate}
	\end{example}

	We observe the following consequence of the domination inequality:
	
	\begin{corollary}\label{homotopy-invariance}
		For any connected Lie group $G$, the $G$--representation volume is a homotopy-type invariant for
		orientable closed smooth manifolds.
	\end{corollary}

\section{Real semisimple Lie groups}\label{Sec-realSemisimpleLieGroups}
	In this section, we recall some standard facts about real semisimple Lie groups.
	See \cite{Ho,He} for general references.
	
	Let $G$ be a connected real Lie group. Denote by $\mathfrak{g}$ the Lie algebra of $G$,
	regarded as the tangent space of $G$ at its neutral element. 
	Denote by $\mathrm{GL}(\mathfrak{g})$ the Lie group of invertible $\R$--linear transformations 
	of $\mathfrak{g}$. The Lie algebra $\mathfrak{gl}(\mathfrak{g})$ of $\mathrm{GL}(\mathfrak{g})$
	can be naturally identified with the Lie algebra of $\R$--endomorphisms $\mathrm{End}(\mathfrak{g})$.
	The adjoint action of $G$ on $\mathfrak{g}$,
	denoted as
		$${\rm Ad}\colon G\to \mathrm{GL}(\mathfrak{g}),$$
	is defined for any $h\in G$ as the derivation (at the neutral element) of the inner automorphism
	$\sigma(h)(g)=hgh^{-1}$, for all $g\in G$. 
	The derivation of ${\rm Ad}$ is hence a Lie algebra homomorphism
		$${\rm ad}\colon \mathfrak{g}\to \mathrm{End}(\mathfrak{g}).$$
	For any $X\in\mathfrak{g}$,  ${\rm ad}(X)$ can be explicitly given by the Lie bracket, 
	namely, ${\rm ad}(X)(Y)=[X,Y]$ for all $Y\in\mathfrak{g}$.
	
	Recall that a Lie algebra $\mathfrak{g}$ is said to be \emph{simple} if it is nonabelian and contains no ideals
	other than $\{0\}$ and $\mathfrak{g}$. It is said to be \emph{semisimple} if it is a direct sum of
	simple Lie algebras. 
	A Lie group $G$ is said to be \emph{semisimple} if its Lie algebra $\mathfrak{g}$ is semisimple.
	
	\subsection{Cohomology}
	Semisimple Lie algebras	can be characterized by cohomology
	through the following \cite[Theorem 24.1]{CE}.
	See Section \ref{Sec-CohomologyOfLieAlgebras} for a review of Lie algebra cohomology.
	
	\begin{theorem}\label{Weyl}
		Let $\mathfrak{g}$ be a finite dimension Lie algebra over any field of characteristic $0$.
		Then $\mathfrak{g}$ is semisimple if and only if $H^*(\mathfrak{g};V)=0$ 
		holds for all irreducible nontrivial $\mathfrak{g}$--modules $V$.
	\end{theorem}
	
	In particular, $H^*(\mathfrak{g};\mathfrak{g}^*)=0$ for real semisimple Lie algebras.
	Moreover, the semisimplicity of $\mathfrak{g}$ implies that any (finite dimensional) $\mathfrak{g}$--module $W$
	is a direct sum of irreducible factors. It follows from Theorem \ref{Weyl}
	that $H^*(\mathfrak{g};W)$ is a direct sum of finitely many (possibly zero) copies of $H^*(\mathfrak{g})$.
	In low dimensions, it can be computed that $H^0(\mathfrak{g})=\R$, and $H^k(\mathfrak{g})=0$
	for $k=1,2,4$. However, there is a canonical nonvanishing class of $H^3(\mathfrak{g})$
	given by the $3$--form $B([\cdot,\cdot],\cdot)$ in $\wedge^3\mathfrak{g}^*$, 
	where $B(X,Y)=\tr({\rm ad}(X){\rm ad}(Y))$ is the Killing form.
	See Chevalley--Eilenberg \cite{CE}.
	
	The following relative version of the vanishing result \cite[Theorem 28.1]{CE} is needed in the rest of this paper. 
	
	\begin{theorem}\label{Weyl-relative}
		Let $\mathfrak{g}$ be a finite dimension Lie algebra over any field of characteristic $0$
		and $\mathfrak{h}$ be a Lie subalgebra.
		If $\mathfrak{g}$ is semisimple, then $H^*(\mathfrak{g},\mathfrak{h};V)=0$ 
		holds for all irreducible nontrivial $\mathfrak{g}$--modules $V$.
	\end{theorem}
	
	\subsection{Center and linearity}
	There are semisimple groups which are not linear, such as $\widetilde{{\rm SL}_2}(\R)$.
	However, it is well known that linearity holds after factoring out the center:
	
	\begin{theorem}\label{semi-simple1}
		Let $Z(G)$ be the center of a connected real semisimple Lie group $G$. Then the following statements hold true:
		\begin{enumerate}
		\item The quotient group $G/Z(G)$ is isomorphic to the identity component of ${\rm Aut}(\mathfrak{g})$ as a Lie group.
		\item The quotient homomorphism $G\to G/Z(G)$ is a covering projection.
		\item The center $Z(G)$ is finitely generated.
		\end{enumerate}
	\end{theorem}	
	
	Here ${\rm Aut}(\mathfrak{g})$ denotes the group of Lie algebra automorphisms of $\mathfrak{g}$, 
	which is contained by the real general linear group ${\rm GL}(\mathfrak{g})$ as a closed (and Zariski closed)
	subgroup. In fact, the isomorphism of the first statement
	can be induced by the adjoint representation ${\rm Ad}\colon G\to \mathrm{GL}(\mathfrak{g})$,
	which is trivial on $Z(G)$. The kernel of ${\rm Ad}$ is exactly $Z(G)$ 
	and the image is the entire identity component of ${\rm Aut}(\mathfrak{g})$ when $G$ is semisimple,
	see \cite[Chapter II, Corollary 5.2 and Corollary 6.5]{He}.
	The second statement is an immediate consequence of the semisimplicity of $G$ 
	that the closed normal subgroup $Z(G)$ has to be $0$--dimensional, and hence discrete.
	It follows that $Z(G)$ is a quotient of the fundamental group $\pi_1(G)$, 
	so the finite generation property of the third statement
	holds by the fact $G$ is homotopy equivalent to any maximal compact subgroup of itself.

	\begin{remark}
	There is a much stronger result asserting that every connected Lie group $G$ admits a universal linear quotient.
	When $G$ is semisimple, the result can be stated as follows.
	Let $\pi\colon \t{G}\to G$ denote the universal covering of $G$ and let $\sigma\colon \t{G}\to G_{\C}$  be the complexification of $\t{G}$. 
	Then every linear representation of $G$ factors through $G/\pi(\ker\sigma)$
	and there exists a faithful linear representation of $G/\pi(\ker\sigma)$,
	which is therefore the universal linear quotient of $G$.
	See \cite[Chapter XVII, Theorem 3.3]{Ho} and \cite{Ho-kernel}.
	\end{remark}
	
	\subsection{Lattices}
	The existence of torsion-free uniform lattices in any connected semisimple Lie group 
	is a celebrated theorem due to A.~Borel \cite[Theorem B]{Borel-lattice},
	which can be paraphrased by the following:
	
	\begin{theorem}\label{Borel}
		Let $G$ be a connected real semisimple Lie group and $H$ be a maximal compact subgroup of $G$.
		Then $G$ has a discrete subgroup $\Gamma$ which acts freely discontinuosly on 
		the homogeneous $G$--space $G/H$ with compact quotient.
	\end{theorem}
	
	In fact, given any torsion-free uniform lattice 
	$\Gamma$ of $G$ by Borel's theorem, the discreteness of $\Gamma$ and the compactness of $H$
	implies that the action of $\Gamma$ on $G/H$ is properly discontinuous.
	Moreover, $\Gamma\cap H$ has to be finite, and therefore trivial as $\Gamma$ contains no torsion.
	The cocompactness of the action is equivalent to the uniformity of $\Gamma$ as $H$ is compact.

\section{Central extension}\label{Sec-centralExtension}
	Let $G$ be a connected semisimple Lie group. 
	We assume for this section that the center $Z$ of $G$ is torsion-free. 
	Denote by $H$ a maximal compact subgroup of $G$ and $X$ the homogeneous space $G/H$.
	In this section, we consider a natural non-splitting central extension $G_{\R}$ of $G$,
	which contains a maximal compact subgroup $H_{\R}$ extending $H$.
	The associated homogeneous space $G_{\R}/H_{\R}$ can be identified with $X$, 
	but there are more transformations coming from the enlarged group $G_{\R}$.
	In our applications, passage from $G$ to $G_{\R}$ allows us more room to deform a representation.

	Note that by our assumption and Theorem \ref{semi-simple1}, the center $Z$ of $G$ is a finitely generated free abelian group.
	Denote by $Z_{\R}$ the (additive) torsion-free abelian Lie group $Z\otimes \R$, which is isomorphic to $\R^{\mathrm{rk}(Z)}$ and contains $Z$
	naturally as the integral lattice. 
	Define the \emph{full central extension} of $G$ to be
		$$G_{\R}=G\times_Z Z_{\R},$$ 
	namely,
	the quotient of $G\times Z_{\R}$ by the $Z$--action 
		$$z\cdot(g,x)\mapsto(gz,x-z\otimes1)$$
	for all $z\in Z$.
	It is clear that $G_{\R}$ naturally contains $G$ as a closed subgroup, which intersects 
	the center $Z_{\R}$ in the naturally embedded $Z$.
	Denote by $\o{G}$ the quotient group $G/Z$. 
	We have the following commutative diagram of Lie group homomorphisms:
	$$\xymatrix{
	\{0\} \ar[r]  & Z \ar[r] \ar[d] & G\ar[r] \ar[d] & \o{G}\ar[r] \ar[d]^{\mathrm{id}} & \{1\} \\
	\{0\} \ar[r]  & Z_{\R} \ar[r]  & G_{\R}	\ar[r] & \o{G}\ar[r]  & \{1\} }$$
	where the rows are exact sequences. 

	\subsection{Representation varieties and lift obstruction}
	For any finitely generated group $\pi$, there are maps between the representation varieties
	$$\xymatrix{
	\c{R}(\pi,G) \ar[r] & \c{R}(\pi,G_{\R}) \ar[r] &\c{R}(\pi,\o{G})
	}$$
	which are naturally induced by the group homomoprhisms
	$$\xymatrix{
	G\ar[r] & G_{\R} \ar[r] &\o{G}.
	}$$
	
	\begin{proposition}\label{rep-var}
	Suppose that $G$ is a connected real semisimple Lie group with torsion-free center $Z$.
	There is a characteristic class for $\o{G}$--representations of finitely generated groups, namely,
	a natural assignment
		$$e_Z\colon \c{R}(\pi,\o{G})\to H^2(\pi;Z)$$
	for any finitely generated group $\pi$.
	Moreover, the following statements are true:
	\begin{enumerate}
	\item The space of representations $\c{R}(\pi,\o{G})$ is a finite union of path-connected components
	of an affine real algebraic variety. 
	The characteristic class $e_Z$ is constant over each path-connected component of $\c{R}(\pi,\o{G})$.
	\item The space of representations $\mathcal{R}(\pi,G)$
		is an $H^1(\pi;Z)$--principal bundle over the union of the path-connected components of
		$\c{R}(\pi,\o{G})$ on which $e_Z$ is trivial.	
	\item The space of representations $\mathcal{R}(\pi,G_{\R})$
		is an $H^1(\pi;Z_{\R})$--principal bundle over the union of the path-connected components of
		$\c{R}(\pi,\o{G})$ on which $e_Z$ is torsion.
	\end{enumerate}	
	\end{proposition}
	
	\begin{proof}
		For any finitely generated group $\pi$ and any $\o{G}$--representation $\eta$ of $\pi$,
		the characteristic class $e_Z(\eta)\in H^2(\pi;Z)$
		can be constructed concretely as follows:
		
		Take any CW complex $K$ which realizes the Eilenberg--MacLane space $K(\pi,1)$.
		Denote by $P_{\eta}$ the associated $\o{G}$--principal bundle $K\times_\eta\o{G}$ over $K$.
		Let $\{U_{\alpha}\}_{\alpha\in A}$ be a locally finite open cover of $K$ over which $P_\eta$ is
		locally trivialized	by sections $s_\alpha\colon U_\alpha\to P_\eta|_{U_\alpha}$.
		The transition functions $g_{\alpha\beta}\in\o{G}$ 
		are given by $s_\beta=s_\alpha g_{\alpha\beta}$ over any $U_\alpha\cap U_\beta$,
		which satisfy the cocycle condition $g_{\alpha\beta}g_{\beta\gamma}g_{\gamma\alpha}=1$
		over any $U_\alpha\cap U_\beta\cap U_\gamma$. Take a lift $\tilde{g}_{\alpha\beta}\in G$
		for each transition function in such a way that $\tilde{g}_{\alpha\beta}=\tilde{g}_{\beta\alpha}^{-1}$.
		It follows that over each $U_\alpha\cap U_\beta\cap U_\gamma$, 
		the lift determines an element
			$$\tilde{g}_{\alpha\beta}\tilde{g}_{\beta\gamma}\tilde{g}_{\gamma\alpha}\in Z.$$
		It is easy to check that this determines a $Z$--valued \v{C}ech $2$--cocycle which represents a cohomology class
			$$e_Z(\eta)\in \check{H}^2(K;Z)\cong H^2(\pi;Z).$$
		Moreover, $e_Z(\eta)$ depends only on $\eta$ and it is functorial 
		with respect homomorphisms between finitely generated groups.
		
		To prove  statement (1), we first observe that $e_Z$ is by definition constant as 
		the representation $\eta$ varies continuously, so it is constant on each path-connected component of 
		$\c{R}(\pi,\o{G})$. To bound the number of the path-connected components,
		we apply Theorem \ref{semi-simple1} about semisimple Lie groups.
		Denote by $\mathfrak{g}$ the Lie algebra of $G$. The adjoint representation
		$\mathrm{Ad}\colon G\to \mathrm{Aut}(\mathfrak{g})$ induces an embedding of $\o{G}$
		as the identity component of $\mathrm{Aut}(\mathfrak{g})$.
		Note that for any continuous path of representations $\rho_t\colon \pi\to \mathrm{Aut}(\mathfrak{g})$,
		the image of every element of $\pi$ must remain in some component of $\mathrm{Aut}(\mathfrak{g})$ all the time.
		It follows that $\c{R}(\pi,\o{G})$ can be identified with those path-connected components
		of $\c{R}(\pi,\mathrm{Aut}(\mathfrak{g}))$ which consist of all the representations into 
		the identity component of $\mathrm{Aut}(\mathfrak{g})$.
		Since $\mathrm{Aut}(\mathfrak{g})$ is a Zariski closed subset of  
		$\mathrm{GL}(\mathfrak{g})$, (thus an affine real algebraic group,)
		$\mathcal{R}(\pi,\mathrm{Aut}(\mathfrak{g}))$ is an affine real algebraic variety
		(namely, a Zariski closed subset of $\R^N$ for some large $N$).
		By H.~Whitney \cite[Theorem 3]{Whitney}, $\c{R}(\pi,\mathrm{Aut}(\mathfrak{g}))$ has at most finitely many
		path-connected components, so the same holds for $\c{R}(\pi,\o{G})$.
		Note that $\o{G}$ is an analytic linear group but not an algebraic group in general.
		
		To prove  statement (2), regard any element $\alpha\in H^1(\pi;Z)$ as a homomorphism $\alpha\colon \pi\to Z$. 
		For any representation $\rho\colon \pi\to G$, we can twist $\rho$
		by $\alpha$, namely, setting $\alpha\cdot\rho\colon \pi\to G$ by $(\alpha\cdot\rho)(g)=\alpha(g)\rho(g)$.
		This induces an action of $H^1(\pi;Z)$ on $\c{R}(\pi,G)$ which is clearly free.
		Note that $Z(G)$ is discrete by Theorem \ref{semi-simple1},
		so the action is discontinuous with respect to the algebraic-convergence topology,
		as can be easily checked on a finite set of generators.
		This action endows $\c{R}(\pi,G)$ with a $H^1(\pi;Z)$--principal bundle structure.
		Every fiber of $\c{R}(\pi,G)$ is clearly projected to a single point of $\c{R}(\pi,\o{G})$.
		As $G/Z$ is connected and covered by $G$, by M.~Culler \cite[Theorem 4.1]{C},
		any continuous path of representations in $\mathcal{R}(\pi,G/Z)$ 
		can be lifted to $\mathcal{R}(\pi,G)$	provided that	some representation of the path can be lifted.
		Therefore, $\c{R}(\pi,G)$ are projected onto a sub-union of path-connected components
		of $\c{R}(\pi,G/Z)$. By the construction of $e_Z$, a $\o{G}$--representation $\eta$ can be lifted
		to a $G$--representation exactly when $e_Z(\eta)\in H^2(\pi;Z)$ vanishes.
		
		Statement (3) can be proved in a similar way. 
		In fact, since $e_Z$ is constant on the components of $\c{R}(\pi,\o{G})$,
		which are at most finitely many, we may take $D$ to be the least common multiple of 
		the orders of those $e_Z$--values which are torsion.
		Using a finite central extension $G_{D^{-1}\Z}=G\times_Z (Z\otimes D^{-1}\Z)$ instead of $G$,
		we see that all the $e_Z$--torsion components of $\c{R}(\pi,\o{G})$ are projected onto
		by $\c{R}(\pi,G_{D^{-1}\Z})$, and hence by $\c{R}(\pi,G_\R)$. 
		Representations of the other components of $\c{R}(\pi,\o{G})$
		cannot be lifted to $\c{R}(\pi,G_\R)$ since $e_Z$ survives in $H^2(\pi;Z_{\R})$.
		The $H^1(\pi;Z_{\R})$--principal bundle structure of $\c{R}(\pi,G_{\R})$
		follows from the same argument as in (2). This completes the proof.
	\end{proof}
	
	%
		%
	%
	
	\subsection{Associated homogeneous space}
	From the proper action of $G$ on its associated contractible homogeneous space $X=G/H$,
	one can extract a concrete model of the central extension $G_\R$ as follows.
		
	\begin{proposition}\label{concreteG_R}
		Let $G$ be a connected semisimple Lie group.  
		Assume that the center $Z$ of $G$ is torsion-free.
		Let $H$ be a maximal compact subgroup of $G$.
		There exists a connected closed torsion-free abelian subgroup $R$ of $G$ with the following properties:
		\begin{itemize}
			\item The subgroup of $G$ generated by $H$ and $R$ is a direct product $HR$.
			\item The intersection of $R$ and $Z$ is a lattice of $R$ and a finite-index subgroup of $Z$.
			\item The center $Z$ of $G$ is contained by $HR$. 
			\item The quotient $HR/Z$ is embedded as a maximal compact subgroup of $\o{G}$
			as induced by the projection.
		\end{itemize}
		Therefore, the center $Z$ of $G$ is embedded as a lattice of $HR/H$ via the quotient projection,
		and there is a uniquely induced isomorphism 
			$$G_{\R}\cong G\times_Z (HR/H).$$
	\end{proposition}
	
	\begin{proof}
		Let $H$ be a maximal compact subgroup of $G$. 
		Let $\o{K}$ be a maximal compact subgroup of $\o{G}$ which contains the projected image of $H$.
		Denote by $K$ the preimage of $\o{K}$ in $G$.
		Since maximal compact subgroups are topologically deformation retracts of any virtually connected Lie group,
		we have homotopy equivalences $H\simeq G$ and $\o{K}\simeq \o{G}$ via the inclusions.
		Therefore, we have $H\simeq K\simeq G$ via the inclusions because $K$ is the covering space of $\o{K}$
		that corresponds the subgroup $\pi_1(G)$ of $\pi_1(\o{G})$.		
		Moreover, $\o{K}$ and $K$ are both connected 
		as $\o{K}$ is homotopy equivalent to the connected group $\o{G}$.
		
		By the structure theory of compact Lie groups, the universal covering group $\t{K}$ of
		the connected compact Lie group $\o{K}$ 
		is isomorphic to a direct product $\t{K}_1\times\cdots\times \t{K}_m\times \R^{n}$ where $\t{K}_i$ are 
		simply-connected simple compact Lie groups, \cite[Theorem 6.6]{HM}.
		The center $Z(\t{K})$ of $\t{K}$
		is the direct product of the finite centers $Z(\t{K}_i)$ and the abelian factor $\R^n$.
		The kernel of the covering projection $\t{K}\to \o{K}$ 
		can be identified with $Z(\t{G})$ by the homotopy equivalence $\o{K}\simeq \o{G}$.
		Denote by $L_K$ the kernel of $\t{K}\to K$.
		There is an induced short exact sequence of discrete abelian groups
		$$\xymatrix{
		\{0\} \ar[r] & L_K \ar[r] & Z(\t{G}) \ar[r] & Z \ar[r] & \{0\},
		}$$
		which splits as $Z$ is assumed to be torsion-free. 
		It follows that the intersection of $L_K$ with $\R^n$ is a finite index subgroup of $L_K$
		which is contained by a unique $\R$--subspace $V_K$ of $\R^n$. 
		Take any lift $\t{Z}$ of $Z$ into $Z(\t{G})$.
		The intersection between $\t{Z}$ and the factor subgroup $\R^n$ 
		is a finite index subgroup of $\t{Z}$, and 
		it is contained by a unique $\R$--subspace $\t{R}$ of $\R^n$.
		However, 
		note that $Z(\t{G})$ may not be a subgroup of the connected component $\R^n$ of $Z(\t{K})$,
		even if it is torsion-free,	so one might not be able to make $\t{Z}$ contained in $\R^n$.
		The fact that $\o{K}$ is compact
		implies the relation
		$$\dim V_K+\dim \t{R}=\dim Z(\t{K}),$$
		or more concretely, there is a direct-sum decomposition $V_K+\t{R}=\R^n$.
		We take the claimed abelian subgroup $R$ of $G$ to be the projected image of $\t{R}$.
		
		Note that the preimage of $H$ in $\t{K}$ is exactly $\t{H}=\t{K}_1\times\cdots\times\t{K}_m\times V_K$,
		which contains the kernel $L_K$ of $\t{K}\to K$.
		Since $\t{R}$ meets $\t{H}$ trivially and centralizes $\t{H}$,
		the double coset $HR$ form a subgroup of $G$ which is the direct product of $H$ and $R$.
		The intersection of $Z$ and $R$ is a finite-index subgroup
		of $Z$ and a lattice of $R$ because the same holds for $\t{Z}$ and $\t{R}$.
		The fact that $Z$ is contained by $HR$ follows from 
		that $\t{Z}$ is contained by $Z(\t{G})\subset Z(\t{H})\t{R}$.
		The quotient $HR/Z$ is a maximal compact subgroup of $\o{G}$ because it is
		exactly $\t{H}\t{R}/Z(\t{G})=\t{K}/Z(\t{G})=\o{K}$.
		Therefore, we have verified the claimed properties about $R$.
		
		To complete the proof, we see that $Z$ is embedded into $HR/H$ as a lattice because it is torsion-free and
		it meets $R$ in a finite index subgroup. The induced isomorphism
			$$G_\R\cong G\times_Z (HR/H)$$
		is immediately implied by the unique isomorphism
			$$Z_{\R}\cong HR/H$$
		which respects to the inclusions of $Z$.
	\end{proof}
	
	In the following, 
	we keep the notation $HR/H$ even though it is canonically isomorphic to $R$.
	We denote the elements of $HR/H$ as cosets $rH$ where $r\in R$.
	In this way, 
	we remind ourselves that $Z$ is not automatically contained by $R$ as between subgroups of $G$.
	
	The homogeneous space $X=G/H$ is equipped with a proper transitive action
		$$G\times_Z(HR/H)\to \mathrm{Diff}(X)$$
	defined by
		$$(g,rH)\cdot(xH)=(gxr^{-1})H.$$
	In this concrete model, there is a distinguished point 
		$$O_X\in X,$$
	namely, the coset $H\in X$.
	The stabilizer of $O_X$, or the isotropy group of the action,
	is the image of the diagonal embedding
		$$\Delta\colon HR/Z\to G\times_Z(HR/H)$$
	defined by $\Delta(hrZ)=(hr,rH)\bmod Z$.
	We denote the isotropy group as	
		$$H_{\R}=\Delta(HR/Z),$$
	while we identify $G_{\R}$ with	$G\times_Z(HR/H)$. 
	Therefore, we have an identification between homogeneous spaces
		$$X\cong G_{\R}/H_{\R}$$
	via the induced equivariant diffeomorphism.
	
	\subsection{Representations and reduction of bundles}\label{Subsec-reduction}
	Let $M$ be a connected closed smooth manifold. Denote by $\t{M}$ the universal cover of $M$.
	For any representation $\rho\co \pi_1(M)\to G_{\R}$,
	after choosing a $\rho$--equivariant smooth developing map $D_\rho\co \t{M}\to X$,
	there is an induced principal $H_{\R}$--bundle $H_{\R}(\t{M},D_\rho)$ over $\t{M}$,
	namely, the pull-back of the canonical principal $H_{\R}$--bundle $G_{\R}\to X$
	via $D_\rho$. 
	The bundle space is the fiber product
		$$H_{\R}(\t{M},D_\rho)=\t{M}\times_{D_\rho} G_{\R},$$
	which consists of those points $(\tilde{m},g)\in \t{M}\times G_{\R}$
	such that $D_\rho(\tilde{m})=g\cdot O_X$.
	The bundle projection takes any point $(\tilde{m},g)$ to $\tilde{m}$.
	Moreover, $H_{\R}(\t{M},D_\rho)$ is equipped with a natural action of $\pi_1(M)$
	given by $\sigma\cdot(\tilde{m},g)\to (\sigma\cdot\tilde{m},\rho(\sigma)\cdot g)$,
	which commutes with the bundle projection,
	so the quotient by the action yields a principal $H_{\R}$--bundle 
	$H_{\R}(M,D_\rho)$ over $M$.
	This is a model of the \emph{$H_{\R}$--reduction} of 
	the flat principal $G_{\R}$--bundle $M\times_\rho G_{\R}$ over $M$,
	induced by the section $s$ of $M\times_\rho X\to M$ that corresponds to $D_\rho$,  
	in the sense that the following diagram of maps commutes:
	$$\xymatrix{
	H_{\R}(M,D_\rho) \ar^{\mathrm{incl.}}[r] \ar[d]& M\times_\rho G_{\R} \ar[d]\\
	M \ar^{s}[r]& M\times_\rho X	
	}$$
	
	The isomorphism type of $H_{\R}$--reductions of flat principal $G_{\R}$--bundles
	$M\times_\rho G_{\R}$	over $M$ depends only on the path-connected component 
	of $\rho\in\mathcal{R}(\pi_1(M),G_{\R})$. However, 
	the models are different as $\rho$ varies. 
	In the following, 
	we exhibit a construction that turns a smooth path of representations 
	into a smoothly varying family of $H_{\R}$--reductions tied to a fixed model.
	The procedure should be routine and easy in principle, 
	but we need to properly formulate and examine a few details.
	
	To keep concrete, we say that a path of representations $\rho_t$ in $\mathcal{R}(\pi_1(M),G_{\R})$ is
	\emph{smooth} if at every element $\sigma\in\pi_1(M)$, 
	the path $\rho_t(\sigma)$ is smooth in $G_{\R}$.
	According to the description of $\mathcal{R}(\pi_1(M),G_{\R})$ from Proposition \ref{rep-var},
	any pair of representations	on a path-connected component of $\mathcal{R}(\pi_1(M),G_{\R})$ can be connected 
	by a piecewise smooth path of representations 
		$$\rho_t\co \pi_1(M)\to G_{\R},$$
	parametrized by $t\in[0,1]$.
		
	\begin{lemma}\label{developing_path}
		Given a smooth path of representations $\rho_t\co \pi_1(M)\to G_{\R}$, parametrized by $t\in[0,1]$, 
		there exists a path of smooth developing maps
			$$D_t\co \t{M}\to X$$
		which is $\rho_t$--equivariant for each $t$ 
		and which varies smoothly with respect to $t$.
	\end{lemma}
	
	\begin{proof}
		Suppose that $\rho_t$ is a smooth path of representations.
		Write $M_{[0,1]}$ for the product $[0,1]\times M$,
		and $M_t$ for any slice $\{t\}\times M$.
		We take the induced smooth bundle 
		over $M_{[0,1]}$ with fibers diffeomorphic to $X$,
		denoted as $M_{[0,1]}\times_\rho X$.
		The total space is the quotient of $[0,1]\times \t{M}\times X$ 
		by the induced action of $\pi_1(M)$, 
		namely, $\sigma\cdot(t,\tilde{m},x)=(t,\sigma\cdot\tilde{m}, \rho_t(\sigma)\cdot x)$.
		The bundle projection is induced by the projection 
		of the universal cover onto the first two factors.
		Since the fiber $X$ is contractible, the bundle admits a continuous section,
		which can be perturbed to be a smooth section
		$s\co M_{[0,1]}\to M_{[0,1]}\times_\rho X$ by standard techniques of approximations
		\cite[Chapter 2, see Section 2, Exercise 3]{Hirsch}.
		Any elevation $\t{s}\co [0,1]\times\t{M}\to [0,1]\times \t{M}\times X$ gives rise to 
		a path of developing maps $D_t\co \t{M}\to X$ such that $D_t(\tilde{m})$
		is the $X$--component of $\t{s}(t,\tilde{m})$. It follows that
		$D_t(\sigma\cdot\tilde{m})=\rho_t(\sigma)\cdot D_t(\tilde{m})$
		for all deck transformations $\sigma\in\pi_1(M)$
		and $D_t$ varies smoothly for $t\in[0,1]$.
		Therefore, $D_t$ is a smooth path of $\rho_t$--equivariant smooth developing maps,
		as claimed.
	\end{proof}
	
	\begin{proposition}\label{model_reduction}
		For every path-connected component $C$ of $\mathcal{R}(\pi_1(M),G_{\R})$,
		there exists a principal $H_{\R}$--bundle 
		$$H_{\R}(M)_C\to M,$$
		whose isomorphism class is determined by $C$.
		Furthermore, 
		with its pull-back to $\t{M}$ denoted by $H_{\R}(\t{M})_C\to \t{M}$,
		the following characterization is true:
		
		Given any smooth path of $\rho_t$--equivariant developing maps
		$D_t\co \t{M}\to X$ with respect to a smooth path of representations $\rho_t$ on $C$,
		parametrized by $t\in[0,1]$,
		there exists a path of $\pi_1(M)$--equivariant isomorphisms of principal $H_\R$--bundles
		$$\phi_t\co H_{\R}(\t{M})_C\to H_{\R}(\t{M},D_t)$$
		which is smooth in the sense that the following composed maps
		vary smoothly with respect to $t$:
		$$H_{\R}(\t{M})_C\stackrel{\phi_t}\longrightarrow H_{\R}(\t{M},D_t)\stackrel{\mathrm{incl}.}\longrightarrow \t{M}\times G_{\R}.$$
	\end{proposition}
	
	\begin{proof}
		Given any smooth path of $\rho_t$--equivariant developing maps
		$D_t\co \t{M}\to X$ with respect to a smooth path of representations $\rho_t$,
		parametrized by $t\in[0,1]$, the principal $H_{\R}$--bundles $H_{\R}(\t{M},D_t)$
		can be put together as a smooth principal $H_{\R}$--bundle
			$$H_{\R}(\t{M}_{[0,1]},D_{[0,1]})\to \t{M}_{[0,1]},$$
		where $\t{M}_{[0,1]}$ is the product $[0,1]\times\t{M}$.
		The bundle space $H_{\R}(\t{M}_{[0,1]},D_{[0,1]})$ consists of the points $(t,\tilde{m},g)$
		of $[0,1]\times \t{M}\times G_{\R}$ with the property that $D_t(\tilde{m})=g\cdot O_X$,
		and the bundle projection takes $(t,\tilde{m},g)$ to $(t,\tilde{m})$. As $\pi_1(M)$ acts freely on
		$H_{\R}(\t{M}_{[0,1]},D_{[0,1]})$ by $\sigma\cdot(t,\tilde{m},g)=(t,\sigma\cdot\tilde{m},\rho_t(\sigma)g)$,
		which commutes with the bundle projection,
		there is an induced smooth principal $H_{\R}$--bundle
			$$H_{\R}(M_{[0,1]},D_{[0,1]})\to M_{[0,1]},$$
		where $M_{[0,1]}$ is the product $[0,1]\times M$.
		There exists a smooth isomorphism of principal $H_{\R}$--bundles
		$$\xymatrix{
		[0,1]\times H_{\R}(M,D_0) \ar^{\Psi}[r] \ar[d]& H_{\R}(M_{[0,1]},D_{[0,1]}) \ar[d]\\
		[0,1]\times M \ar^{\mathrm{id}}[r]& M_{[0,1]}		
		}$$
		which can be constructed as follows.
		
		Observe that both $[0,1]\times H_{\R}(M,D_0)$ and $H_{\R}(M_{[0,1]},D_{[0,1]})$
		restricts to be $H_{\R}(M,D_0)$ over the $0$--slice $M_0$ of the base space.
		As $M_0$ is a deformation retract of $M_{[0,1]}$, there is a continuous isomorphism $\Psi$
		of principal $H_{\R}$--bundles which fits into the claimed commutative diagram.
		We can modify $\Psi$ to be smooth	by the following argument:		
		Cover $M$ by compact disks $U_1,\cdots,U_N$, such that each $U_i$
		is contained by an open regular neighborhood $U'_i$. The product of $U_i$ and $U'_i$
		with $[0,1]$ are denoted as $U_{i,[0,1]}$ and $U'_{i,[0,1]}$ accordingly.
		Then the smooth principal $H_{\R}$--bundles $H_{\R}(\t{M}_{[0,1]},D_{[0,1]})$
		and $[0,1]\times H_{\R}(\t{M},D_0)$
		can be trivialized over $U'_{i,[0,1]}$ 
		by smooth local sections $s_{i,[0,1]}$ and $\mathrm{id}_{[0,1]}\times s_{i,0}$
		accordingly, such that the $H_{\R}$--valued transition functions on overlaps 
		are all smooth for both bundles.
		With the local trivializations given by the local sections, the continuous bundle-isomorphism 
		$\Psi$ is determined by its local functions
		$h_i\co U_{i,[0,1]}\to H_{\R}$.
		Those functions are defined by the relation
		$$\Psi\left((\mathrm{id}_{[0,1]}\times s_{i,0})(t,\tilde{m})\right)=s_{i,[0,1]}(t,\tilde{m})\cdot h_i(t,\tilde{m}),$$
		where the right action of $h_i(t,\tilde{m})$ applies only to the $G_{\R}$--component
		by right multiplication.
		Using the standard approximation techniques \cite[Chapter 2]{Hirsch},
		we can slightly perturb $h_1$ over $U'_{1,[0,1]}$ to make it smooth in an open neighborhood of $U_1$.
		Modify $\Psi$ accordingly over $U'_{1,[0,1]}$.
		Proceeding inductively, we can modify the next $h_i$ over $U'_{i,[0,1]}$ 
		without affecting the already modified $\Psi$ on $U_1\cup\cdots\cup U_{i-1}$.
		After finitely many steps, we eventually obtain a modified bundle-isomorphism $\Psi$
		which is smooth.
		
		In particular, we see that the restriction of $\Psi$ to the $t$--slice $\{t\}\times H_{\R}(M,D_0)$ induces 
		a smooth path of smooth isomorphisms of principal $H_{\R}$--bundles over $M$:
			$$\Psi_t\co H_{\R}(M,D_0)\stackrel{\cong}\longrightarrow H_{\R}(M,D_t).$$
		After passing to covering spaces, 
		it gives rise to a smooth path of smooth $\pi_1(M)$--equivariant isomorphisms of principal $H_{\R}$--bundles over $\t{M}$:
			$$\psi_t\co H_{\R}(\t{M},D_0)\cong H_{\R}(\t{M},D_t).$$
		
		To conclude the proof, it suffices to show that
		given any path-connected component $C$ of $\mathcal{R}(\pi_1(M),G_{\R})$, 
		any representation $\rho_C$ and a smooth $\rho_C$--equivariant developing map
		$D_C\co \t{M}\to X$, the isomorphism class of the principal $H_{\R}$--bundle
		$H_{\R}(M,D_C)$ depends only on $C$.
		Taking this fact for granted for the moment, we can fix any such model
		as the claimed $H_{\R}(M)_C$.
		In fact, given any smooth path of $\rho_t$--equivariant developing maps
		$D_t\co \t{M}\to X$ with respect to a smooth path of representations $\rho_t$ on $C$,
		parametrized by $t\in[0,1]$, the claimed smooth path of $\pi_1(M)$--equivariant isomorphisms
		$\phi_t$ can be taken as the following path of composed isomorphisms
			$$H_{\R}(\t{M})_C\stackrel{\cong}\longrightarrow H_{\R}(\t{M},D_0)\stackrel{\psi_t}\longrightarrow H_{\R}(\t{M},D_t),$$
		where the first isomorphism does not depend on $t$.
		
		To see that the isomorphism class of $H_{\R}(M,D_C)$ depends only on $C$,
		we argue as the following.
		For any fixed $\rho_C$ of $C$, different choices of smooth $\rho_C$--equivariant developing maps
		can be joined by a smooth path of $\rho_C$--equivariant developing maps, for example,
		using the argument of Lemma \ref{choices} and the developing map interpretation of sections.
		So the isomorphism type of $H_{\R}(M,D_C)$ depends only on $\rho_C$.
		For any other choice $\rho'_C$ of $C$, and any $D'_C$ accordingly,
		there is a piecewise smooth path connecting
		$\rho_C$ and $\rho'_C$. So by Lemma \ref{developing_path} applied to each smooth piece,
		and the induced isomorphism $\psi_1$ constructed above, we obtain an isomorphism between
		$H_{\R}(M,D_C)$ and $H_{\R}(M,D'_C)$, through a finite sequence of intermediate isomorphisms.
		Therefore, the isomorphism type of $H_{\R}(M,D_C)$ depends only on $C$.
		This completes the proof.		
	\end{proof}

\section{Cohomology of Lie algebras}\label{Sec-CohomologyOfLieAlgebras}
	Cohomology of Lie algebra supplies a useful device for computation of volume for representations.
	From the level of forms it can be derived 
	through several standard operations assoicated with the exterior algebra.
	The constructions that we recall in this section apparently work for arbitrary ground fields,
	but we specialize our discussion to real Lie algebras.
	
	\subsection{Operations on forms}
	Let $\mathfrak{g}$ be a real Lie algebra. The exterior algebra
		$$A^*(\mathfrak{g})=\bigoplus_k\wedge^k\mathfrak{g}^*$$
	is a graded algebra over $\R$. We regard  any exterior $k$--form over $\mathfrak{g}$
	as an antisymmetric $\R$--multilinear function of $k$ variables in $\mathfrak{g}$ and valued in $\R$.
	More generally, given any $\mathfrak{g}$--module $V$, namely,
	a $\R$--vector space with a Lie algebra homomorphism $\mathfrak{g}\to\mathrm{End}(V)$,
	we can consider $V$--valued exterior $k$--forms over $\mathfrak{g}$.
	Such twisted forms constitute 
	$$A^*(\mathfrak{g};V)=\bigoplus_k{\rm Hom}(\wedge^k\mathfrak{g},V),$$
	which is a graded left $A^*(\mathfrak{g})$--module with respect to the wedge product.
	Elements of ${\rm Hom}(\wedge^k\mathfrak{g},V)$ are regarded as antisymmetric $\R$--multilinear maps
	$f\co\mathfrak{g}\times\cdots\times\mathfrak{g}\to V$. 
		
	On any $A^*(\mathfrak{g};V)$ there are three natural operations, 
	namely, the differential, Lie derivatives, and the interior product.
	The differential $\delta\co A^k(\mathfrak{g};V)\to A^{k+1}(\mathfrak{g};V)$ is defined by
	\begin{eqnarray*}
	(\delta f)(X_1,\cdots,X_{k+1})&=&\sum_{i}(-1)^{i+1}X_i\cdot f(X_1,\cdots,\hat{X_i},\cdots,X_{k+1})+\\
	&&\sum_{i<j}(-1)^{i+j}f([X_i,X_j],X_1,\cdots,\hat{X_i},\cdots,\hat{X_j},\cdots,X_{k+1}).
	\end{eqnarray*}
	The differential of $A^*(\mathfrak{g})$ is denoted particularly as $\ud$,
	and there the first summation is gone.
	Given any $X\in\mathfrak{g}$, 
	the Lie derivative $L_X\co A^k(\mathfrak{g};V)\to A^k(\mathfrak{g};V)$ 
	is defined by 
	$$(L_X f)(X_1,\cdots,X_k)= X\cdot f(X_1,\cdots,X_k)+\sum_i f(X_1,\cdots,[X_i,X],\cdots,X_k),$$ 
	and 
	the interior product $i_X\co A^k(\mathfrak{g};V)\to A^{k-1}(\mathfrak{g};V)$ is defined by 
	$$(i_X f)(X_1,\cdots,X_{k-1})=f(X,X_1,\cdots,X_{k-1}).$$
	Recall that $L_X,i_X$ and $\delta$ are related by 
	\begin{equation}\label{LR}
	L_X=\delta\circ i_X+i_X\circ \delta.
	\end{equation}
	
	If we take a basis $e_1,\cdots,e_n$ of $\mathfrak{g}$ and 
	let $e^*_1,\cdots,e^*_n$ denote the dual basis of $\mathfrak{g}^*$,
	then the usual differential can be calculated by 
	\begin{equation}\label{delta'}
	2\ud\omega=\sum_i e^*_i\wedge L_{e_i}\omega,
	\end{equation}
	for $\omega\in A^*(\mathfrak{g})$, \cite[Ch.~I, Sec.~1.1, Formula (7)]{BW}.
	In particular, if $\omega$ is a $1$--form, then   
  \begin{equation}\label{delta1} 
	\ud\omega=\sum_{i<j}\omega([e_j, e_i])\,e^*_i\wedge e^*_j.
	\end{equation}
	Here and throughout this paper, we adopt the wedge product notation
	$$e^*_{i_1}\wedge\cdots\wedge e^*_{i_k}=\sum_{\sigma\in\mathfrak{S}_k}{\rm sgn}(\sigma)\,e^*_{i_{\sigma(1)}}\otimes\cdots\otimes e^*_{i_{\sigma(k)}}.$$

	\subsection{Relative cohomology with module coefficients}
	With the notations above, it can be checked that
	$(A^*(\mathfrak{g};V),\delta)$ becomes a cochain complex
	over $\R$, and the induced cohomology is denoted by 
	$$H^*(\mathfrak{g};V)=H^*(A^*(\mathfrak{g};V)).$$
	When $V$ is trivially $\R$, the cohomology is usually called the cohomology
	of the Lie algebra $\mathfrak{g}$, denoted simply by $H^*(\mathfrak{g})$.
	
	More generally, suppose that $\mathfrak{h}$ is a Lie subalgebra of $\mathfrak{g}$.
	We denote by 
	$$A^*(\mathfrak{g},\mathfrak{h};V)\subset A^*(\mathfrak{g};V)$$
	the forms that are annihilated by both $L_X$ and $i_X$ for every $X\in\mathfrak{h}$.
	Note that $A^*(\mathfrak{g},\mathfrak{h};V)$
	is closed under $\delta$ by the formula (\ref{LR}), 
	so $\delta$ restricts to be a differential on $A^*(\mathfrak{g},\mathfrak{h};V)$. 
	This allows us to define the relative cohomology:
	$$H^*(\mathfrak{g},\mathfrak{h};V)=H^*(A^*(\mathfrak{g},\mathfrak{h};V)).$$
		
	If $(\mathfrak{g},\mathfrak{h})$ are the Lie algebras of a real Lie group $G$ and a connected closed subgroup $H$,
	the relative cohomology $H^*(\mathfrak{g},\mathfrak{h})$ can be identified with the cohomology of $G$--invariant differential forms
	of the homogeneous space $X=G/H$. Namely, denote by 
	$A^*(X)^G$ the graded algebra of $G$--invariant differential forms of $X$, and by $H^*(A^*(X)^G)$ its corresponding cohomology.
	The evaluation at the origin yields an isomorphism $A^*(X)^G\to A^*(\mathfrak{g},\mathfrak{h})$ that commutes with the differential,
	so it induces an identification	
	$$H^*(A^*(X)^G)\cong H^*(\mathfrak{g},\mathfrak{h}),$$
  see \cite[Chapter 1, Section 1.6]{BW}.
	
	\subsection{The coadjoint module}\label{Subsec-coadjoint}
	The dual vector space $\mathfrak{g}^*$ can be endowed with a coadjoint representation 
	$$-{\rm ad}^*\co\mathfrak{g}\to{\rm End}(\mathfrak{g}^*)$$
	to become a $\mathfrak{g}$--module.
	The action $\mathfrak{g}\times\mathfrak{g}^*\to \mathfrak{g}^*$ is therefore defined by 
	\begin{equation}\label{contra} 
	(X,\omega)\mapsto L_X\omega=\omega([\cdot,X]).
	\end{equation}
		
	\begin{remark}\label{contragradient}
		Some authors adopt a different convention, using the coaction ${\rm ad}^*$ of $\mathfrak{g}$ on $\mathfrak{g}^*$.
		The difference is that $-\mathrm{ad}^*\co \mathfrak{g}\to \mathrm{End}(\mathfrak{g}^*)$ is a Lie algebra homomorphism which
		preserves the Lie bracket. In fact, the Jacobi identity implies that
		\begin{eqnarray}\label{NLH}
			{\rm ad}^*([X,Y])={\rm ad}^*(Y)\circ{\rm ad}^*(X)-{\rm ad}^*(X)\circ{\rm ad}^*(Y).
		\end{eqnarray}
		For this reason, our formula of the differential $\delta$ is not the same as \cite[(23.1), page 115]{CE},
		but the convention here agrees with \cite[Chapter 1]{BW}.
	\end{remark}
	%
	
	There exists a natural homomorphism of graded $\R$--modules of degree $-1$:	
	$$J\co A^*(\mathfrak{g})\to A^{*-1}(\mathfrak{g};\mathfrak{g}^*)$$
	 defined by the formula:
	\begin{equation}\label{Jdef}
	(J\omega)(X_1,\cdots,X_{k-1})(X)=i_X(\omega(\cdot,X_1,\cdots,X_{k-1}))=\omega(X,X_1,\cdots,X_{k-1}),
	\end{equation}
	for any $\omega\in A^k(\mathfrak{g})$.	
	
	\begin{lemma}\label{JRelation}
	The homomorphism $J$ satisfies the following commutative relations:
	\begin{eqnarray}
		\delta\circ J&=&-J\circ \ud\\
		i_X\circ J&=&-J\circ i_X\\
		L_X\circ J&=&J \circ L_X
	\end{eqnarray}	
	\end{lemma}
	
	\begin{proof}
		Assuming the first two relations, we can derive the third relation by:
		$$L_X\circ J=(\delta\circ i_X+i_X\circ\delta)\circ J=J\circ (\ud\circ i_X+i_X\circ\ud)= J\circ L_X.$$
		For any $\omega\in A^k(\mathfrak{g})$, the first two relations can be verified directly as follows.
		For any $X_1,\cdots,X_k,Y\in\mathfrak{g}$, we compute:
		\begin{eqnarray*}
		&&\left((\delta J\omega)(X_1,\cdots,X_k)\right)(Y)\\
		&=&\sum_{i}(-1)^{i+1}\left(-\mathrm{ad}^*(X_i)\cdot((J\omega)(X_1,\cdots,\widehat{X_i},\cdots,X_k))\right)(Y)+\\
		&&\sum_{i<j}(-1)^{i+j}\left((J\omega)([X_i,X_j],X_1,\cdots,\widehat{X_i},\cdots,\widehat{X_j},\cdots,X_k)\right)(Y)\\
		&=&\sum_{i}(-1)^{i+1}\omega([Y,X_i],X_1,\cdots,\widehat{X_i},\cdots,X_k)+\\
		&&\sum_{i<j}(-1)^{i+j}\omega(Y,[X_i,X_j],X_1,\cdots,\widehat{X_i},\cdots,\widehat{X_j},\cdots,X_k)\\
		&=&(-1)\times\sum_{i}(-1)^{1+(i+1)}\omega([Y,X_i],X_1,\cdots,\widehat{X_i},\cdots,X_k)+\\
		&&(-1)\times\sum_{i<j}(-1)^{(i+1)+(j+1)}\omega([X_i,X_j],Y,X_1,\cdots,\widehat{X_i},\cdots,\widehat{X_j},\cdots,X_k)\\		
		&=&-(\ud\omega)(Y,X_1,\cdots,X_k)\\
		&=&-\left((J\ud\omega)(X_1,\cdots,X_k)\right)(Y).
		\end{eqnarray*}
		This shows $\delta\circ J=-J\circ\ud$.
		For any $X,X_1,\cdots,X_{k-1},Y\in\mathfrak{g}$, we compute:
		\begin{eqnarray*}
		\left((i_X J\omega)(X_1,\cdots,X_{k-1})\right)(Y)
		&=&\left((J\omega)(X,X_1,\cdots,X_{k-1})\right)(Y)\\
		&=&\omega(Y,X,X_1,\cdots,X_{k-1}),
		\end{eqnarray*}
		and
		\begin{eqnarray*}
		\left((J\,i_X\omega)(X_1,\cdots,X_{k-1})\right)(Y)
		&=&(i_X\omega)(Y,X_1,\cdots,X_{k-1})\\
		&=&\omega(X,Y,X_1,\cdots,X_{k-1}).
		\end{eqnarray*}		
		This shows $i_X\circ J=-J\circ i_X$.
	\end{proof}
	
	As an immediate consequence,
	for any Lie subalgebra $\mathfrak{h}$ of $\mathfrak{g}$,
	there are an induced homomorphism of differential graded $\R$--modules of degree $-1$:
	$$J\co A^*(\mathfrak{g},\mathfrak{h})\to A^{*-1}(\mathfrak{g},\mathfrak{h};\mathfrak{g}^*),$$
	and an induced homomorphism of graded $\R$--modules of degree $-1$:
	$$J'\co H^*(\mathfrak{g},\mathfrak{h})\to H^{*-1}(\mathfrak{g},\mathfrak{h};\mathfrak{g}^*).$$

\section{Volume rigidity of representations}\label{Sec-volumeRigidityOfRepresentations}
	In this section, we prove that for the full central extension of a connected real semisimple Lie group
	with torsion-free center, the volume of representations for any given manifold is locally constant.
		
	\begin{theorem}\label{rigidity}
		Let $G$ be a connected real semisimple Lie group. 
		Suppose that the center $Z$ of $G$ is torsion-free,
		and denote by $G_{\R}$ the full central extension of $G$. 
		For any closed oriented smooth manifold $M$, the volume of representations
			$$\mathrm{vol}_{G_\R}\colon \mathcal{R}(\pi_1(M),G_{\R})\to \R$$
		is constant on every path-connected component of 
		the representation space $\mathcal{R}(\pi_1(M),G_{\R})$.		
	\end{theorem}
	
	The rest of this section is devoted to the proof of Theorem \ref{rigidity}.
	Fix a maximal compact subgroup $H$ of $G$ and denote by $X=G/H$
	the associated contractible homogeneous $G$--space.
	We identify the proper action of $G_\R$ on $X$,
	following Proposition \ref{concreteG_R}, so that the extended isotropy group is
	a maximal compact subgroup $H_\R$ of $G_\R$.
	Fix a $G_{\R}$--invariant volume form	$\omega_X$ of $X$.
	Therefore, the volume of representations of this section are considered with respect to the triple
		$$(G_{\R},X,\omega_X)$$
	in accordance with Notation \ref{usualNotation}.
	Suppose that $M$ is a closed oriented manifold of the same dimension as $X$,
	as in other dimensions the volume is constantly zero by our convention.
	
	Let $C$ be any path-connected component of the representation space $\mathcal{R}(\pi_1(M),G_{\R})$.
	According to the description of $\mathcal{R}(\pi_1(M),G_{\R})$ from Proposition \ref{rep-var},
	any pair of representations	on $C$ can be connected 
	by a piecewise smooth path of representations. 
	Therefore, it suffices to show that $\mathrm{vol}_{G_{\R}}$ is constant restricted to
	any smooth path of representations on $C$.
		
	\subsection{Derivative of developing paths}\label{Subsec-derivative_of_developing}
	Suppose that 
		$$\rho_t\co \pi_1(M)\to G_{\R}$$
	is a smooth path of representations on $C$, parametrized by $t\in[0,1]$.
	By Lemma \ref{developing_path}, there exists 
	a smooth path of smooth $\rho_t$--equivariant developing maps
		$$D_t\co \t{M}\to X.$$
	It follows that there is an induced smooth path of homomorphisms between differential graded $\R$--algebras:
		$$D^\sharp_t\co A^*(X)^{G_{\R}}\to A^*(\t{M})^{\pi_1(M)},$$
	namely, $\R$--linear homomorphisms which preserves the wedge product, the exterior differential,
	and the grading.
	The $G_{\R}$--invariant differential forms on $X$ are naturally identified with 
	exterior forms on the tangent space $T_{H}X$ at the coset $H$ of 
	the homogeneous space $X=G/H$, where the isotropy group is $H_{\R}$.
	The $\pi_1(M)$--invariant differential forms on $\t{M}$ are naturally identified with
	the pull-backs of the differential forms on $M$ via the covering.
	Therefore, the induced smooth path of homomorphisms between differential graded $\R$--algebras
	can be identified as
		$$D^\sharp_t\co A^*(\mathfrak{g}_{\R},\mathfrak{h}_{\R})\to A^*(M).$$
	The derivative of $D^\sharp_t$ with respect to $t$ gives rise to
	a smooth path of homomorphisms between differential graded $\R$--modules
		$$\dot{D}^\sharp_t\co A^*(\mathfrak{g}_{\R},\mathfrak{h}_{\R})\to A^*(M).$$
	On any open subset $U$ of $M$ with local coordinates $(u_1,\cdots,u_m)$,
		$$D^\sharp_t(\omega)=\sum_{i_1<\cdots<i_k}f^\omega_{i_1,\cdots,i_k}(t,u_1,\cdots,u_m)\,\ud u_{i_1}\wedge\cdots\wedge \ud u_{i_k},$$
	where the coefficient functions are linear in $\omega\in A^k(\mathfrak{g}_{\R},\mathfrak{h}_{\R})$
	and smooth in $(t,u_1,\cdots,u_m)$,
	and
		$$\dot{D}^\sharp_t(\omega)=\sum_{i_1<\cdots<i_k}
		\left.\frac{\partial f^\omega_{i_1,\cdots,i_k}}{\partial t}\right|_{(t,u_1,\cdots,u_m)}\,\ud u_{i_1}\wedge\cdots\wedge \ud u_{i_k}.$$
	It follows that there are induced homomorphisms of graded $\R$--algebras
		$$D^*_t\co H^*(\mathfrak{g}_{\R},\mathfrak{h}_{\R})\to H^*(M;\R),$$
	and homomorphisms of graded $\R$--modules
		$$\dot{D}^*_t\co H^*(\mathfrak{g}_{\R},\mathfrak{h}_{\R})\to H^*(M;\R).$$
		
	\begin{lemma}\label{isConstant}
		Given a class $[\omega]\in H^*(\mathfrak{g}_{\R},\mathfrak{h}_{\R})$, 
		the smooth path of classes
		$D^*_t[\omega]$ is constant in $H^*(M;\R)$ if and only if 
		the derivative classes $\dot{D}^*_t[\omega]$ vanishes for all $t$.
	\end{lemma}
	
	\begin{proof}
		It is obvious from the local expression that derivation and integration with respect to $t$
		commute with the differentials of the differential graded $\R$--modules
		$A^*(\mathfrak{g}_{\R},\mathfrak{h}_{\R})$ and $A^*(M)$. As it can be checked completely locally,
		the conclusion is a simple fact of multivariable calculus.
	\end{proof}
	
	\subsection{Factorization of derivatives}\label{Subsec-factorization}
	
	\begin{lemma}\label{factorization}
		Given any smooth path of representations $\rho_t\in\mathcal{R}(\pi_1(M),G_{\R})$,
		parametrized by $t\in[0,1]$,
		and any smooth path of smooth $\rho_t$--equivariant developing maps
		$D_t\co \t{M}\to X$,
		there exists a factorization of $\dot{D}^\sharp_t$
		as a composition of homomorphisms, of degree $-1$ and $+1$,
		between differential graded $\R$--modules
		$$A^*(\mathfrak{g}_{\R},\mathfrak{h}_{\R})\stackrel{J}{\longrightarrow}
		A^{*-1}(\mathfrak{g}_{\R},\mathfrak{h}_{\R};\mathfrak{g}^*_{\R})\stackrel{\dot{F}_t}{\longrightarrow}
		A^*(M).$$
		Therefore, there exists a factorization of $\dot{D}^*_t$ accordingly
		$$H^*(\mathfrak{g}_{\R},\mathfrak{h}_{\R})\stackrel{J'}{\longrightarrow}
		H^{*-1}(\mathfrak{g}_{\R},\mathfrak{h}_{\R};\mathfrak{g}^*_{\R})\stackrel{\dot{F}'_t}{\longrightarrow}
		H^*(M).$$ 
	\end{lemma}
	
	\begin{proof}
		Denote by $C$ the path-connected component of $\mathcal{R}(\pi_1(M),G_{\R})$ which contains the path $\rho_t$.
		Denote by $H_{\R}(M,D_t)\to M$ the principal $H_{\R}$--bundles which are the $H_{\R}$--reductions
		induced by $D_t$, of the flat $G_{\R}$--bundles $M\times_{\rho_t}X\to M$ with fiber $X$,
		(Subsection \ref{Subsec-reduction}).
		
		Fix a model of the $H_{\R}$--reduction associated with $C$,
		namely, a principal $H_{\R}$--bundle $H_{\R}(M)_C\to M$
		as asserted by Proposition \ref{model_reduction}.
		With the pull-backs to $\t{M}$ denoted by $H_{\R}(\t{M})_C$ and $H_{\R}(M,D_t)$
		accordingly, there are smooth isomorphisms of $\pi_1(M)$--equivariant principal $H_{\R}$--bundles
		$\phi_t$, as guaranteed by Proposition \ref{model_reduction}, which fits into the following
		sequence of smooth morphisms of $\pi_1(M)$--equivariant principal $H_{\R}$--bundles:
		$$\xymatrix{
		H_{\R}(\t{M})_C \ar^{\phi_t}[r] \ar[d] & H_{\R}(\t{M},D_t) \ar^{\mathrm{incl.}}[r] \ar[d] &
		\t{M}\times G_{\R} \ar^{\mathrm{proj.}}[r] \ar[d] & G_{\R} \ar[d] \\
		\t{M} \ar^{\mathrm{id}}[r] & \t{M} \ar^{\mathrm{id}\times D_t}[r] & 
		\t{M}\times X \ar^{\mathrm{proj.}}[r] & X  
		}$$
		(The action of $\pi_1(M)$ on any $G_{\R}$ or $X$ factors are understood as by $\rho_t$.)
		The composition along the lower row of the commutative diagram is nothing but the developing maps $D_t$.
		We denote by
			$$\Delta_t\co H_{\R}(\t{M})_C\to G_{\R}$$
		the maps obtained by composition along the upper row,
		which are therefore $\rho_t$--equivariant, 
		and vary smoothly as parametrized by $t\in[0,1]$.
		The induced maps of invariant differential forms
			$$\Delta_t^\sharp\co A^*(G_{\R})^{G_{\R}}\to A^*(H_{\R}(\t{M})_C)^{\pi_1(M)}$$
		can be identified with the homomorphisms between differential graded $\R$--algebras
			$$\Delta_t^\sharp\co A^*(\mathfrak{g}_{\R})\to A^*(H_{\R}(M)_C).$$
		The derivative of $\Delta_t^\sharp$ with respect to $t$ are
		homomorphisms between differential graded $\R$--modules
			$$\dot{\Delta}_t^\sharp\co A^*(\mathfrak{g}_{\R})\to A^*(H_{\R}(M)_C).$$
		Via the bundle projection of the principal $H_{\R}$--bundles,
		the derivatives $\dot{D}_t^\sharp$ can be identified with 
		the restriction of $\dot{\Delta}_t^\sharp$, as seen by the commutative diagram
		of homomorphisms between differential graded $\R$--modules
		$$\xymatrix{
		A^*(\mathfrak{g}_{\R}) \ar^{\dot{\Delta}_t^\sharp}[r] & A^*(H_{\R}(M)_C)\\
		A^*(\mathfrak{g}_{\R},\mathfrak{h}_{\R}) \ar^{\dot{D}_t^\sharp}[r] \ar[u]& A^*(M) \ar[u]
		}$$
		where the vertical homomorphisms are inclusions,
		and forms of $A^*(M)$ are identified with forms of $A^*(H_{\R}(M)_C)$ that are
		$H_{\R}$--horizontal and $H_{\R}$--invariant.
		
		Note that $H_{\R}(M)_C$ is a principal $H_{\R}$--bundle.
		For any $X\in\mathfrak{h}_{\R}$,
		the flow on $H_{\R}(M)_C$ given by the right action of 
		the $1$--parameter subgroup $\exp(sX)$ gives rise to a vector field of $H_{\R}(M)_C$,
		still denoted as $X$ for simplicity.
		Then the interior product $i_X$ and Lie derivative $L_X$ are operators on $A^*(H_{\R}(M)_C)$.
		Therefore, a form of $A^*(H_{\R}(M)_C)$ is $H_{\R}$--horizontal if and only if it is annihilated by $i_X$
		for all $X\in\mathfrak{h}_{\R}$, and $H_{\R}$--invariant if and only if it is annhilated by $L_X$.
		
		There is a natural factorization of $\dot{\Delta}_t^\sharp$ 
		into homomorphisms between differential graded $\R$--modules
		$$A^*(\mathfrak{g}_{\R}) \stackrel{J}\longrightarrow A^{*-1}(\mathfrak{g}_{\R};\mathfrak{g}_{\R}^*)
		\stackrel{F_t}\longrightarrow A^*(H_{\R}(M)_C),$$
		which is just the Leibniz rule of derivatives.
		More precisely, the homomorphism $J$ is the canonical homomorphism
		introduced in Subsection \ref{Subsec-coadjoint}, (Formula \ref{Jdef}).
		Identifying $A^{k-1}(\mathfrak{g}_{\R};\mathfrak{g}_{\R}^*)$ with $(\wedge^{k-1}\mathfrak{g}_{\R}^*)\otimes\mathfrak{g}_{\R}^*$,
		the homomorphism $F_t$ is defined by specifying
			$$F_t(\xi\otimes\omega)=\dot{\Delta}_t^\sharp(\omega)\wedge \Delta_t^\sharp(\xi).$$
		The factorization is readily checked by 
		\begin{eqnarray*}
			\dot{\Delta}_t^\sharp(\omega_1\wedge\cdots\wedge\omega_k)&=&
			\sum_{i=1}^{k}\Delta_t^\sharp(\omega_1)\wedge\cdots\wedge \dot{\Delta}_t^\sharp(\omega_i)\wedge\cdots\wedge \Delta_t^\sharp(\omega_k)\\
			&=&
			\sum_{i=1}^{k}(-1)^{i+1}\dot{\Delta}_t^{\sharp}(\omega_i)\wedge
			\Delta_t^\sharp\left(\omega_1\wedge\cdots\wedge \widehat{\omega_i}\wedge\cdots\wedge \omega_k\right)\\
			&=&F_t(J(\omega_1\wedge\cdots\wedge\omega_k)).
		\end{eqnarray*}
		Furthermore, it is easy to verify that 
		$$\ud\circ F_t=-F_t\circ \delta,$$
		and that for any $X\in\mathfrak{h}$,
		$$i_X\circ F_t=-F_t\circ i_X,$$
		and 
		$$L_X\circ F_t=F_t\circ L_X.$$
		This means that $F_t$ is indeed a homomorphism of differential graded $\R$--modules, 
		and its restriction induces a homomorphism of differential graded $\R$--modules
		$$F_t\co A^{*-1}(\mathfrak{g}_{\R},\mathfrak{h}_{\R};\mathfrak{g}_{\R}^*)\to A^*(M).$$
		The restriction of $J$ is a homomorphism of differential graded $\R$--modules
		$$J\co A^{*}(\mathfrak{g}_{\R},\mathfrak{h}_{\R})\to A^*(\mathfrak{g}_{\R},\mathfrak{h}_{\R};\mathfrak{g}_{\R}^*),$$
		(Subsection \ref{Subsec-coadjoint}).
		Passing to the restrictions, therefore, we obtain a factorization
		of the derivative $\dot{D}_t^\sharp$ into homomorphism of differential graded $\R$--modules
		$$A^*(\mathfrak{g}_{\R},\mathfrak{h}_{\R}) \stackrel{J}\longrightarrow A^{*-1}(\mathfrak{g}_{\R},\mathfrak{h}_{\R};\mathfrak{g}_{\R}^*)
		\stackrel{F_t}\longrightarrow A^*(M),$$
		as claimed.
		The cohomological factorization is an immediate consequence of the factorization on the form level.
		This completes the proof.
	\end{proof}
	
	\subsection{A vanishing lemma}
	
	\begin{lemma}\label{vanishing}
		$$H^{\dim(X)-1}(\mathfrak{g}_{\R},\mathfrak{h}_{\R};\mathfrak{g}_{\R}^*)=0.$$
	\end{lemma}
	
	\begin{proof}
		We compute the cohomology $H^1(\mathfrak{g}_{\R},\mathfrak{h}_{\R};\mathfrak{g}_{\R})$
		directly and derive the claimed vanishing by Poincar\'{e} duality.
		Here the coefficient $\mathfrak{g}_{\R}$ is considered as the $\mathfrak{g}_{\R}$--module
		with the adjoint action.
	
		First observe that
			$$H^1(\mathfrak{g},\mathfrak{h})=0.$$
		In fact, for any $\omega\in A^1(\mathfrak{g},\mathfrak{h})$ we have
		$(\ud\omega)(X,Y)=-\omega([X,Y])$. If $\ud\omega=0$, the semisimplicity of $\mathfrak{g}$
		implies $[\mathfrak{g},\mathfrak{g}]=\mathfrak{g}$, so $\omega=0$.
		Using the decomposition of $\mathfrak{g}_{\R}$ into the center and the simple factors,
			$$\mathfrak{g}_{\R}=\mathfrak{z}(\mathfrak{g}_{\R})\oplus \mathfrak{g}_1\oplus\cdots\oplus \mathfrak{g}_r,$$		
		we obtain an induced decomposition of $H^1(\mathfrak{g},\mathfrak{h};\mathfrak{g}_{\R})$
		into the direct sum of $H^1(\mathfrak{g},\mathfrak{h})\otimes_{\R}\mathfrak{z}(\mathfrak{g}_{\R})$
		and $H^1(\mathfrak{g},\mathfrak{h};\mathfrak{g}_i)$, 
		which are all trivial by the above observation and Theorem \ref{Weyl-relative}.
		Therefore,
			$$H^1(\mathfrak{g},\mathfrak{h};\mathfrak{g}_{\R})=0.$$
				
		The cohomology classes of $H^1(\mathfrak{g}_{\R},\mathfrak{h}_{\R};\mathfrak{g}_{\R})$
		are by definitions represented by $\R$--linear maps
			$$f\co \mathfrak{g}_{\R}\to \mathfrak{g}_{\R}$$
		with the property that 
			$$i_X(f)=f(X)=0$$
		for all $X\in\mathfrak{h}_{\R}$ and
			$$(\delta f)(X,Y)=[X,f(Y)]-[Y,f(X)]-f([X,Y])=0$$ 
		for all $X,Y\in\mathfrak{g}_{\R}$ and
			$$L_X(f)(Y)=[X,f(Y)]-f([X,Y])=0$$
		for all $X\in\mathfrak{h}_{\R}$ and $Y\in\mathfrak{g}_{\R}$.
		For any such $f$, the restriction of $f$ to $\mathfrak{g}$
		is a $1$--cocycle of $A^1(\mathfrak{g},\mathfrak{h};\mathfrak{g}_{\R})$.
		
		When $\mathfrak{h}$ is nontrivial,
		the restriction of $f$ to $\mathfrak{g}$ is $0$ because
		$A^0(\mathfrak{g},\mathfrak{h};\mathfrak{g}_{\R})$ is trivial
		and $H^1(\mathfrak{g},\mathfrak{h};\mathfrak{g}_{\R})=0$.
		Then $f$ is trivial on $\mathfrak{g}_{\R}$ since it is also trivial on $\mathfrak{h}_{\R}$.
		
		When $\mathfrak{h}$ is trivial, the vanishing cohomology $H^1(\mathfrak{g};\mathfrak{g}_{\R})=0$
		implies that there exists some $0$--cochain $U\in \mathfrak{g}_{\R}$, and
			$$f(X)=[X,U]$$
		for all $X\in \mathfrak{g}$.
		We argue that $f(X)=0$ if $X$ lies in the center $\mathfrak{z}(\mathfrak{g}_{\R})$.
		In fact, the formula of $\delta f$ implies that $[Y,f(X)]=0$ for all $Y\in\mathfrak{g}_{\R}$
		if $X\in\mathfrak{z}(\mathfrak{g}_{\R})$. 
		So $f(X)\in \mathfrak{z}(\mathfrak{g}_{\R})$.
		On the other hand, recall that	$\mathfrak{h}_{\R}$ is 
		the diagonal of the subalgebra $\mathfrak{r}\oplus\mathfrak{z}(\mathfrak{g}_\R)$,
		where $\mathfrak{r}$ is the Lie algebra of the subgroup $R$ of $G$,
		according to the description of $H_{\R}$ from Proposition \ref{concreteG_R}.
		It follows that for any $X\in\mathfrak{z}(\mathfrak{g}_{\R})$,
		there exists some $X'\in\mathfrak{r}$ such that $X+X'\in\mathfrak{h}_{\R}$.
		So we have $f(X)=-f(X')=-[X',U]$, which lies in $\mathfrak{g}$.
		It follows that $f(X)$ lies in $\mathfrak{z}(\mathfrak{g}_{\R})\cap\mathfrak{g}=\{0\}$
		whenever $X$ lies in $\mathfrak{z}(\mathfrak{g}_{\R})$.
		This shows 
			$$f(X)=[X,U]=(\delta U)(X)$$
		for all $X\in\mathfrak{g}_{\R}$. In other words, any $1$--cocycle $f$ is a coboundary
		$\delta U$ of some $0$--cochain $U$.
		
		Therefore, $H^1(\mathfrak{g}_{\R},\mathfrak{h}_{\R};\mathfrak{g}_{\R})$ always vanishes.
		By the Poincar\'e duality \cite[Chapter I, Proposition 1.5]{BW},
		$$H^{\dim(X)-1}(\mathfrak{g}_{\R},\mathfrak{h}_{\R};\mathfrak{g}_{\R}^*)\cong
		H^1(\mathfrak{g}_{\R},\mathfrak{h}_{\R};\mathfrak{g}_{\R})^*=0,$$
		which completes the proof.		
	\end{proof}
	
	\subsection{Proof of Theorem \ref{rigidity}}
	We summarize the discussion so far and finish the proof of Theorem \ref{rigidity}.	
	Let $G$ be a connected semisimple Lie group with torsion-free center $Z$.
	Adopt the notations $(G_{\R},X,\omega_X)$ as before. 
	For any closed oriented smooth manifold $M$ of the dimension the same as $X$,
	and for any path-connected component $C$
	of $\mathcal{R}(\pi_1(M),G_{\R})$, it suffices to show that
	$\mathrm{vol}_{G_{\R}}(M,\rho_t)$ is constant for every smooth path of representations
	on $C$, as we have argued at the beginning of this section.
	
	Take a smooth path of $\rho_t$--equivariant developing maps $D_t\co \t{M}\to X$
	as in Subsection \ref{Subsec-derivative_of_developing}. Adopting the notations there,
	we have
	$$\mathrm{vol}_{G_{\R}}(M,\rho_t)=\int_{\mathcal{F}} D_t^*\omega_X=\langle D_t^*[\omega_X], [M]\rangle.$$
	Here the integral is over a fixed fundamental domain $\mathcal{F}$ of $\t{M}$,
	and the pairing is the canonical pairing 
	between the top-dimensional cohomology and homology of $M$.
	The derivative homomorphism
	$$\dot{D}_t^*\co H^{\dim(X)}(\mathfrak{g}_{\R},\mathfrak{h}_{\R};\R)\to H^{\dim(X)}(M;\R)$$ 
	is constantly zero as it factors through $H^{\dim(X)-1}(\mathfrak{g}_{\R},\mathfrak{h}_{\R};\mathfrak{g}^*_{\R})=0$,
	by Lemmas \ref{factorization} and \ref{vanishing}.
	It follows that the classes $D_t^*[\omega_X]$ in $H^{\dim(X)}(M;\R)$ are constant
	independent of $t$, (Lemma \ref{isConstant}).
	Therefore, the volume $\mathrm{vol}_{G_{\R}}(M,\rho_t)$ is constant 
	as $t$ varies.
	
	This completes the proof of Theorem \ref{rigidity}.

\section{Cocompact coextension}\label{Sec-cocompactCoextension}
	In this section, we study the structure of connected real Lie groups that contain
	closed cocompact connected semisimple Lie subgroups.
	
	\begin{proposition}\label{cocompact_coextension}
		Let $G$ be a connected real Lie group. 
		Suppose that $G$ contains a closed cocompact connected semisimple Lie subgroup,
		then there exists an exact sequence of homomorphisms of Lie groups
		$$\{0\}\longrightarrow Z(G)_{\mathtt{tor}}\longrightarrow G\longrightarrow \hat{G}_{\R}\longrightarrow T\longrightarrow \{0\}$$
		where $Z(G)_{\mathtt{tor}}$ is the maximal compact central subgroup of $G$, 
 $T$ is a connected compact abelian Lie group, and $\hat{G}_{\R}$
		is the full central extension of a connected semisimple Lie group $\hat{G}$ with torsion-free center.
	\end{proposition}
	
	In fact, the Lie groups $\hat{G}_{\R}$, and $T$ are determined by $G$ up to isomorphism, 
	as the proof indicates. 
	Note also that the maximal compact central subgroup $Z(G)_{\mathtt{tor}}$
	is a possibly disconnected closed normal subgroup of $G$.
	Before we prove Proposition \ref{cocompact_coextension},
	we point out a partial converse as the following:
	
	\begin{lemma}\label{cocompact_coextension_if}
		Suppose that $G$ is 
		a real connected Lie group that fits into an exact sequence as of Proposition \ref{cocompact_coextension}.
		If $Z(G)_{\mathtt{tor}}$ is finite,
		then $G$ contains a closed cocompact connected semisimple normal subgroup.
	\end{lemma}
	
	\begin{proof}
		Since $\hat{G}$ is a closed cocompact connected semisimple normal subgroup of $\hat{G}_{\R}$,
		which necessarily lies in the image of $G$ by the exact sequence,
		the identity component of the preimage of $\hat{G}$ in $G$ yields 
		a closed cocompact connected semisimple normal subgroup of $G$ as claimed.
	\end{proof}
	
	The rest of this section is devoted to the proof of Proposition \ref{cocompact_coextension}.
	To this end, suppose that $G_1$ is a closed connected semisimple Lie subgroup of a real Lie group $G$, 
	such that the coset space $G/G_1$ is compact.
	
	We first show that $G$ has a reductive Lie algebra, and that the commutator subgroup $[G,G]$ is
	a connected closed cocompact semisimple Lie subgroup of $G$. 
	Then we derive the exact sequence by studying the universal covering $\t{G}\to G$.

	By the structure theory of Lie groups, the universal covering group $\t{G}$ of $G$ is a Lie-group 
	semidirect product of its maximal connected solvable subgroup 
	with a maximal simply-connected semisimple subgroup, as induced by the Lie-algebra decomposition
	as a semidirect product, (see \cite[Chapter XII, Theorem 1.2, and Chapter XI]{Ho}).
	The kernel of the covering is a discrete central subgroup of $\t{G}$ \cite[Chapter I, Exercise 1]{Ho}.		
	Therefore, the maximal normal solvable subgroup $G_{\mathtt{sol}}$ of $G$ is closed, and possibly disconnected,
	and the quotient group $\o{G}_{\mathtt{ss}}$ is a connected semisimple real Lie group with trivial center.
	
	Denote by $\o{G}_1$ the image of $G_1$ in $\o{G}_{\mathtt{ss}}$. The preimage of $\o{G}_1$ in $G$
	is obviously the double-coset $G_{\mathtt{sol}}G_1$, which forms a subgroup of $G$ that contains $G_{\mathtt{sol}}$.	
	We have a commutative diagram of Lie-group homomorphisms
	$$\xymatrix{
	\{1\} \ar[r] &G_{\mathtt{sol}} \ar[r]  & G \ar[r] & \o{G}_{\mathtt{ss}} \ar[r]& \{1\}\\
	\{1\} \ar[r] &G_{\mathtt{sol}} \ar[r] \ar[u]^{\mathrm{id}} & G_{\mathtt{sol}}G_1 \ar[r] \ar[u]& \o{G}_1 \ar[r]\ar[u]& \{1\}
	}$$
	where the horizontal arrows have closed image and the rows are short exact sequences.
	The following lemma implies that vertical arrows are all closed maps as well.
	
	\begin{lemma}\label{oG1}
		The subgroup $\o{G}_1$ of $\o{G}_{\mathtt{ss}}$ is closed and cocompact.
	\end{lemma}
	
	\begin{proof}
		The intersection of $G_1\cap G_{\mathtt{sol}}$ is a closed normal solvable
		subgroup of $G_1$. It is necessarily discrete and central in $G_1$,
		because the Lie algebra of that subgroup
		has to be the trivial ideal by the semisimplicity of $G_1$.
		It follows that $\o{G}_1$	is a connected (analytic) semisimple subgroup of $\o{G}_{\mathtt{ss}}$.
		Since $\o{G}_{\mathtt{ss}}$ is semisimple with trivial center,
		it is a closed real algebraic linear group (Theorem \ref{semi-simple1}).
		In this case, it is known that $\o{G}_1$ must be a closed subgroup of $\o{G}_{\mathtt{ss}}$,
		\cite[Chapter II, Exercise and Further Results D.4]{He}.
	\end{proof}
		
	\begin{lemma}\label{Gsol}
		The maximal solvable normal subgroup $G_{\mathtt{sol}}$ of $G$ is abelian. 
		In other words, the Lie algebra of $G$ is reductive.
	\end{lemma}
	
	\begin{proof}	
		By Lemma \ref{oG1}, the subgroup $G_{\mathtt{sol}}G_1$, 
		which equals the preimage of $\o{G}_1$ in $G$, 
		is a closed subgroup of $G$.
		By the assumption, $G_1$ is a closed subgroup of $G$, 
		so the projection of $G$ onto the coset space $G/G_1$
		is an analytic open mapping, \cite[Chapter II Section 4]{He}. 
		It follows that the image $(G_{\mathtt{sol}}G_1)/G_1$ is closed in $G/G_1$.
		By the assumption, $G/G_1$ is compact, so $(G_{\mathtt{sol}}G_1)/G_1$ 
		is a compact (possibly disconnected) manifold without boundary.
		The component of $(G_{\mathtt{sol}}G_1)/G_1$ that contains the identity coset $G_1$
		can be identified with the coset space $G^\circ_{\mathtt{sol}}/(G^\circ_{\mathtt{sol}}\cap G_1)$,
		where $G^\circ_{\mathtt{sol}}$	the identity component of $G_{\mathtt{sol}}$.
		
		Observe that $G^\circ_{\mathtt{sol}}\cap G_1$ is virtually central in $G^\circ_{\mathtt{sol}}$,
		or in other words, a finite-index subgroup of $G^\circ_{\mathtt{sol}}\cap G_1$
		is central in $G^\circ_{\mathtt{sol}}$.
		In fact, the adjoint action induces a linear representation $\alpha\co G_1\to\mathrm{Aut}(\mathfrak{g}_{\mathtt{sol}})$,
		and any element of $G^\circ_{\mathtt{sol}}\cap G_1$ is central in
		$G^\circ_{\mathtt{sol}}$ if and only if it lies in the kernel of $\alpha$.		
		However, as $G_1$ is semisimple, 
		the kernel of $\alpha$ meets the center of $G_1$ in a finite-index subgroup of that center,
		by \cite[Chapter XVII Theorem 3.3 and Theorem 2.1]{Ho}.
		Since $G^\circ_{\mathtt{sol}}\cap G_1$ is central in $G_1$,
		some finite-index subgroup of $G^\circ_{\mathtt{sol}}\cap G_1$ is central in $G^\circ_{\mathtt{sol}}$.
		
		Therefore, $G^\circ_{\mathtt{sol}}\cap G_1$ is a cocompact discrete virtually central subgroup of $G^\circ_{\mathtt{sol}}$,
		and the coset space $G^\circ_{\mathtt{sol}}/(G^\circ_{\mathtt{sol}}\cap G_1)$ is a quotient Lie group.
		The universal cover $\t{G^\circ_{\mathtt{sol}}}$ is a simply-connected solvable Lie group, which is in particular
		contractible. It follows from the structure of compact Lie groups that 
		$G^\circ_{\mathtt{sol}}/(G^\circ_{\mathtt{sol}}\cap G_1)$ is abelian.
		In particular, 
		the Lie algebra of the maximal solvable normal subgroup $G_{\mathtt{sol}}$ of $G$ is abelian,
		so the Lie algebra of $G$ is reductive.
	\end{proof}
	
	By Lemma \ref{Gsol}, the commutator subgroup of $G$ is the semisimple part of $G$,
	namely, it is the maximal semisimple analytic subgroup of $G$.
	We denote
		$$G_{\mathtt{ss}}=[G,G].$$
	Note that $G_{\mathtt{ss}}$ covers $\o{G}_{\mathtt{ss}}$ under the projection of $G$.
	Denote by $K$ the maximal compact connected normal subgroup of $G_{\mathtt{ss}}$,
	and by $\o{K}$ its projected image in $\o{G}_{\mathtt{ss}}$.
	We also observe that $G_1$ is contained by $G_{\mathtt{ss}}$,
	since the corresponding Lie algebras satisfy
		$$\mathfrak{g}_1=[\mathfrak{g}_1,\mathfrak{g}_1]\leq[\mathfrak{g},\mathfrak{g}]=\mathfrak{g}_{\mathtt{ss}}.$$
	The double coset $KG_1$ forms a subgroup of $G$ as $K$ is normal.
	
	\begin{lemma}\label{KG1}
		The commutator subgroup $G_{\mathtt{ss}}$ of $G$ equals $KG_1$. 
		Therefore, $G_{\mathtt{ss}}$ is closed and cocompact in $G$.
	\end{lemma}
	
	\begin{proof}
		We observe that the subgroup $KG_1$ is closed and cocompact in $G$. In fact, since $G_1$ is closed and $K$ is compact,
		$(KG_1)/G_1$ is compact in the coset space $G/G_1$, so $KG_1$ is closed in $G$. 
		Since $G_1$ is cocompact in $G$, the subgroup $KG_1$ is cocompact as well.
		
		Therefore, it suffices to show that the Lie algebra of $KG_1$ equals the Lie algebra of $G_{\mathtt{ss}}$,
		or equivalently, that $\o{K}\o{G}_1$ equals $\o{G}_{\mathtt{ss}}$.
		
		Because $\o{G}_1$ is closed in $\o{G}_{\mathtt{ss}}$ (Lemma \ref{oG1}),
		and $\o{K}$ is compact in $\o{G}_{\mathtt{ss}}$, the coset space $(\o{K}\o{G}_1)/\o{G}_1$
		is compact in $\o{G}_{\mathtt{ss}}/\o{G}_1$.
		It follows that $\o{K}\o{G}_1$ is a closed and cocompact subgroup of $\o{G}_{\mathtt{ss}}$.
		Since $\o{K}$ is normal and compact, it suffices to show
		that the quotient Lie group $(\o{K}\o{G}_1)/\o{K}$ equals $\o{G}_{\mathtt{ss}}/\o{K}$.
		
		In other words, possibly after passing to the quotient,
		we may assume that the semisimple Lie group $\o{G}_{\mathtt{ss}}$ contains no compact normal subgroups.
		We must show under this assumption that 
		any connected cocompact closed semisimple subgroup $\o{G}_1$ equals $\o{G}_{\mathtt{ss}}$.
		This is a consequence of Borel's Density Theorem.
		In fact, consider the adjoint representation
		$\mathrm{Ad}\co \o{G}_{\mathtt{ss}}\to\mathrm{Aut}(\mathfrak{g}_{\mathtt{ss}})$.
		Since the subspace $\mathfrak{g}_1$ is invariant under the action of $\o{G}_1$,
		Borel's Density Theorem \cite[Theorem 5.28]{Raghunathan} implies that $\mathfrak{g}_1$ 
		is invariant under the action of $\o{G}_{\mathtt{ss}}$. This means that $\mathfrak{g}_1$ 
		is an ideal of $\mathfrak{g}_{\mathtt{ss}}$, or equivalently,
		$\o{G}_1$ is a normal subgroup of $\o{G}_{\mathtt{ss}}$.
		The quotient Lie group $\o{G}_{\mathtt{ss}}/\o{G}_1$ has to be trivial, since otherwise
		it would be a semisimple Lie group of noncompact type, which is impossible provided that $\o{G}_1$ is
		cocompact in $\o{G}_{\mathtt{ss}}$.
		This shows that $\o{G}_1$ equals $\o{G}_{\mathtt{ss}}$, which completes the proof.
	\end{proof}
		
	Therefore, we see that $G$ has a reductive Lie algebra, and its commutator subgroup 
	is closed and cocompact. 
			
	\begin{lemma}\label{killing_central_torsion}
		The quotient Lie group $G/Z(G)_{\mathtt{tor}}$ has a reductive Lie algebra,
		and	its commutator subgroup is closed and compact.
		Furthermore, the center of $G/Z(G)_{\mathtt{tor}}$ is torsion-free.
	\end{lemma}
	
	\begin{proof}
		The center of $G/Z(G)_{\mathtt{tor}}$ is obviously torsion-free, otherwise
		$Z(G)_{\mathtt{tor}}$ would not be the maximal compact central subgroup.
		The Lie algebra of $G/Z(G)_{\mathtt{tor}}$ is reductive as the quotient only factors out
		some abelian ideal of the Lie algebra of $G$.
		The preimage of the commutator subgroup of $G/Z(G)_{\mathtt{tor}}$ in $G$ is 
		the double coset $Z(G)_{\mathtt{tor}}G_{\mathtt{ss}}$.
		Note that $Z(G)_{\mathtt{tor}}G_{\mathtt{ss}}$ forms a subgroup of $G$
		since $Z(G)_{\mathtt{tor}}$ is normal. 
		Moreover, it is the preimage of the image of $Z(G)_{\mathtt{tor}}$ in the coset space 
		$G/G_{\mathtt{ss}}$, so the cocompactness of $G_{\mathtt{ss}}$ implies 
		that $Z(G)_{\mathtt{tor}}G_{\mathtt{ss}}$ is closed in $G$. 
		It is cocompact in $G$ since $G_{\mathtt{ss}}$ is already cocompact. 
		Therefore, the image of $(Z(G)_{\mathtt{tor}}G_{\mathtt{ss}})/Z(G)_{\mathtt{tor}}$ in the quotient group
		$G/Z(G)_{\mathtt{tor}}$ is a closed and cocompact subgroup as asserted.
	\end{proof}
	
	With Lemma \ref{killing_central_torsion}, we may argue for $G/Z(G)_{\mathtt{tor}}$
	instead of $G$, so the claimed exact sequence can be obtained by the following lemma.
	
	\begin{lemma}\label{hatGR}
		Suppose in addtion that the center of $G$ is torsion-free.
		Then there exists a homomorphism of Lie groups
			$$G\to \hat{G}_{\R}$$
		where $\hat{G}_{\R}$ is the full central extension 
		of a connected real semisimple Lie group $\hat{G}$ with torsion-free center.
		Moreover, the homomorphism is injective and the image is a closed and cocompact normal subgroup
		of $\hat{G}_{\R}$ that contains $\hat{G}$.
	\end{lemma}
	
	\begin{proof}
		Denote by $\t{G}$ the universal covering Lie group of $G$, and by $\Lambda$ the kernel of the covering projection,
		which is a discrete, central subgroup of $\t{G}$.
		Denote by $G'$ the commutator subgroup of $G$, which is closed and cocompact in $G$.
		Since the Lie algebra of $G$ is reductive, $\t{G}$ is the Lie-group direct product
		of the identity component of the center $Z^\circ(\t{G})$ and the universal covering Lie group $\t{G}'$ of $G'$.
		Note that $\t{G}'$ can be identified with the commutator subgroup of $\t{G}$,
		and it is the maximal connected semisimple subgroup of $\t{G}$,
		which is a direct product of simply-connected simple Lie subgroups.
		
		By the assumption that $G$ has no central torsion,
		the identity component of the center $Z^\circ(\t{G})$ is isomorphic to $\R^{\dim(Z^\circ(\t{G}))}$.
		Moreover, $\Lambda$ contains the torsion subgroup $Z(\t{G}')_{\mathtt{tor}}$ of the center $Z(\t{G}')$ of $\t{G}'$.
		Note that $Z(\t{G}')$ is a finitely generated abelian group.
		We take
			$$\t{G}'_{\mathtt{free}}=\t{G}'/Z(\t{G}')_{\mathtt{tor}},$$
		namely, the characteristic free abelian covering group of $G'$.
		Note that the center $Z(\t{G}'_{\mathtt{free}})$ is discrete and torsion-free.
		It follows from the above description that the quotient of $\t{G}$ by $Z(\t{G}')_{\mathtt{tor}}$
		is a direct product $\t{G}'_{\mathtt{free}}\times Z^{\circ}(\t{G})$. 
		
		The kernel of the covering projection $\t{G}'_{\mathtt{free}}\times Z^{\circ}(\t{G})\to G$ 
		is the free central subgroup, 
			$$\Lambda_{\mathtt{free}}\subset Z(\t{G}'_{\mathtt{free}})\times Z^{\circ}(\t{G}),$$
		namely, the quotient group $\Lambda/Z(\t{G}')_{\mathtt{tor}}$.
		Again by the assumption that $G$ has no central torsion,
		the quotient of $Z(\t{G}'_{\mathtt{free}})\times Z^{\circ}(\t{G})$ by $\Lambda_{\mathtt{free}}$
		is torsion-free abelian.
		This implies that the projection of $\Lambda_{\mathtt{free}}$ to $Z(\t{G}'_{\mathtt{free}})$ 
		is a direct summand of $Z(\t{G}'_{\mathtt{free}})$.
		The assumption that $G'$ is closed and cocompact in $G$ 
		implies that the projection of $\Lambda_{\mathtt{free}}$
		to $Z^{\circ}(\t{G})$ is discrete and cocompact in $Z^{\circ}(\t{G})$,
		or in other words, the projection is a lattice of $Z^\circ(\t{G})$.
		This implies a direct-sum decomposition 
			$$Z(\t{G}'_{\mathtt{free}})=A\oplus B\oplus C$$
		such that $A\oplus B$ is isomorphic with $\Lambda_{\mathtt{free}}$ via the projection $\Lambda_{\mathtt{free}}\to Z(\t{G}'_{\mathtt{free}})$
		and that $B$ is the identified with the kernel of the projection $\Lambda_{\mathtt{free}}\to Z^{\circ}(\t{G})$.
		Note that $B$ is contained by $\Lambda_{\mathtt{free}}$, and the rank of $A$ equals $\dim(Z^\circ(\t{G}))$.
		We take a connected semisimple Lie group
			$$\hat{G}=\t{G}'_{\mathtt{free}}/B,$$
		which covers $G'$ and has torsion-free center isomorphically projected to $A\oplus C$.		
		The above description shows that there exists an isomorphism of Lie groups
			$$G\cong \hat{G}\times_A (A\otimes\R).$$
		The right-hand side can be canonically embedded into $\hat{G}\times_{A\oplus C}((A\oplus C)\otimes \R)$,
		which is the full central extension $\hat{G}_{\R}$ of $\hat{G}$.
		It is clear that this embedding has closed and cocompact image.
		It is normal because it already contains $\hat{G}$.
		Therefore, the Lie group $\hat{G}$ is as desired.	
	\end{proof}
	
	By Lemmas \ref{killing_central_torsion} and \ref{hatGR},
	we obtain an exact sequence
	$$\{1\}\longrightarrow G/Z(G)_{\mathtt{tor}}\longrightarrow \hat{G}_{\R}\longrightarrow T\longrightarrow \{0\},$$
	where $T$ is the connected compact abelian Lie group which is the cokernel
	of the homomorphism in the middle.
	Combining with the canonical short exact sequence
	$$\{0\}\longrightarrow Z(G)_{\mathtt{tor}}\longrightarrow G\longrightarrow G/Z(G)_{\mathtt{tor}}\longrightarrow \{1\},$$
	we obtain the exact sequence as claimed.
	This completes the proof of Proposition \ref{cocompact_coextension}.

\section{Volume finiteness and nontriviality}\label{Sec-volumeFinitenessAndNontriviality}
	In this section, we prove Theorem \ref{main-volume},
	namely, the finiteness and nontriviality of the volume function	
	for real connected Lie groups that contains a closed cocompact semisimple subgroup.
	
%
	
	\begin{proof}[{Proof of Theorem \ref{main-volume}}]
		Since $G$ contains a closed and cocompact semisimple Lie subgroup $G_1$
		we can take a torsion-free uniform lattice $\Gamma$ of $G_1$, 
		by Borel's theorem on the existence of uniform lattices (Theorem \ref{Borel}).
		Then $\Gamma$ is also uniform in $G$. 
		Take $M$ to be the orbit space of $X$ quotient by $\Gamma$.
		Note that $M$ is aspherical as it is covered by a contractible space $X$.
		The inclusion $\rho_0\co \Gamma\to G$ is 
		a discrete faithful representation of $\pi_1(M)\cong\Gamma$,
		and
			$$\mathrm{vol}_{G}(M,\rho_0)=\int_{M}\omega_X>0.$$
		Therefore, the volume function $\mathrm{vol}_{G}$ is nontrivial 
		for some closed oriented smooth manifold $M$.
		
		It remains to show that $\mathrm{vol}_{G}$ takes only finitely many values for any $M$.
		By formula for the volume of induced representations (Proposition \ref{change_of_volume} (2))
		and the characterization of $G$ (Proposition \ref{cocompact_coextension}),
		it suffices to show that $\mathrm{vol}_{G_{\R}}$
		takes only finitely many values for the full central extension $G_{\R}$
		of any real connected semisimple Lie group $G$ with torsion-free center.
		
		In this case, by Proposition \ref{rep-var}, there are at most finitely many path-connected components
		of the space of representations $\mathcal{R}(\pi_1(M);G_{\R})$.
		By Theorem \ref{rigidity}, the volume function is constant on every path-connected component.
		Therefore, the volume function $\mathrm{vol}_{G_{\R}}$ 
		takes only finitely many values.
		This completes the proof.
	\end{proof}

\section{Conclusions}\label{Sec-conclusions}
In conclusion, for any connected Lie group that contains a closed cocompact semisimple subgroup,
the representation volume of that group is a finite nontrivial homotopy-type invariant
for closed orientable smooth manifolds. In this case, the representation volume is virtually essentially
the representation volume associated to the semisimple subgroup of the considered Lie group.
For semisimple Lie groups with infinite center, the associated representation volumes
are useful invariants to detect Property D for certain manifolds with vanishing simplicial volume.

In the following we propose a few further questions.

\begin{question}
	In every dimension $n\geq 4$, is there an orientable closed aspherical smooth $n$--manifold $M$
	with Property D, whereas the representation volume $\mathrm{V}(M',G)$ 
	vanishes for every possible Lie group $G$ and every finite cover $M'$ of $M$? 
	Here we require that $G$ defines a finite representation volume 
	$V(-,G)$ in the dimension $n$.
\end{question}

For dimension $3$, any orientable closed $3$--manifold with property D has a finite cover
with nonvanishing Seifert volume, see \cite{DLSW}.
In dimension 4, the similar question for the simplicial volume is open \cite{FL}.
On the other hand, surface bundles over surfaces are known to have the simplicial volume 
$\|X\|$ at least $4\chi(X)=4\chi(F)\chi(B)>0$, provided that $\chi(F)<0$ and $\chi(B)<0$
for the fiber $F$ and the base $B$, \cite{HK}. 
However, very little is known about their virtual representations.

\begin{question}
	Are there simply-connected closed manifolds of any dimension that has Property D?
\end{question}

Highly connected manifolds tend to fail Property D.
For example, it is known that any $(n-1)$--connected $2n$--manifold $M$ admits self-maps
of nonzero degree, \cite{DuW}. In particular, 
admissible self-maping degrees of any simply-connected $4$--manifold
are all the perfect squares $0,1,4,9,16,\cdots$.
For simply-connected manifolds, representation volumes are by definition zero 
so they provide no obstruction.

\begin{question}
	What is the largest class of connected Lie groups that define nontrivial representation volumes?
\end{question}

\end{document}